\documentclass[12pt]{article}

\usepackage{fullpage,amsfonts,amsmath,amsthm,amssymb,mathrsfs,graphicx,upgreek}
\usepackage[normalem]{ulem}
\usepackage{verbatim,url}
\usepackage{bbm}
\usepackage{graphicx,hyperref,enumitem}
\usepackage[noadjust,sort]{cite}
\usepackage{tocloft}

\usepackage[margin=1in]{geometry}

\makeatletter
\def\blfootnote{\xdef\@thefnmark{}\@footnotetext}
\makeatother

\usepackage{tikz}
\usetikzlibrary{positioning}

\newcommand{\be}{\begin{equation}}
	\newcommand{\ee}{\end{equation}}

\newtheorem{theorem}{Theorem}[section]

\newtheorem{Conjecture}[theorem]{Conjecture}
\newtheorem{lemma}[theorem]{Lemma}
\newtheorem*{Definition}{Definition}

\newtheorem{claim}[theorem]{Claim}
\theoremstyle{definition}
\newtheorem{Remark}[theorem]{Remark}

\def\Time{\mathcal{T}}

\newcommand{\WZtable}[1]{\vspace{-1ex}\[
\begin{tabular}{c|c|c}
  condition & $W^*$ & $\mathcal{Z}$ \\
  \hline\vphantom{$\tilde W$}#1\end{tabular}\]}

\def\G{{\mathcal G}}
\def\po{\operatorname{\mathbf{Po}}}
\def\Bin{\operatorname{\mathbf{Bin}}}

\def\Pr{\operatorname{\mathbb P}}
\def\range#1{\operatorname{rng}(#1)}

\newcommand{\norm}[1]{{\left\|#1\right\|}}
\def\bfd{{\dvec}}
\def\time{{\tau}}

\def\low{{G_{*}}}
\def\up{{G^{*}}}

\def\eps{\varepsilon}
\def\abs#1{\lvert#1\rvert} 
 
\def\norm#1{\lVert#1\rVert}

\def\dfrac#1#2{\lower0.15ex\hbox{\large$\textstyle\frac{#1}{#2}$}}
\def\({\bigl(}
\def\){\bigr)}
\def\st{\mathrel{:}}

\def\nicebreak{\vskip 0pt plus 50pt\penalty-300\vskip 0pt plus -50pt }

\def\S{\boldsymbol{S}}

\def\calG{\mathcal{G}}

\def\degvec{\boldsymbol{\operatorname{deg}}}

\def\dvec{\boldsymbol{d}}
\def\evec{\boldsymbol{e}}
\def\vvec{\boldsymbol{v}}
\def\tvec{\boldsymbol{t}}

\def\svec{\boldsymbol{s}}
\def\rvec{\boldsymbol{r}}

\def\xvec{\boldsymbol{x}}
\def\yvec{\boldsymbol{y}}

\def\gvec{\boldsymbol{g}}

\def\trans{^{\mathrm{T}}\!}

\def\E{\operatorname{\mathbb{E}}}

\def\Var{\operatorname{Var}}

\def\Reals{{\mathbb{R}}}

\def\Naturals{{\mathbb{N}}}

\let\leq=\leqslant
\let\ge=\geqslant
\let\geq=\geqslant

\numberwithin{equation}{section}

\date{}

\title{Kim--Vu's sandwich conjecture is true for  
	$d \gg \log^4 n$
	\blfootnote{
		This paper combines and corrects two preliminary versions: one appeared in Proc SODA 2020 and arxiv:1906.02886, and the other appeared  in  arxiv:2011.09449.}
}
\author{
	Pu Gao\thanks{Research{\vrule height 2.5ex width 0ex} supported by ARC DE170100716, ARC DP160100835, and NSERC RGPIN-04173-2019.} \\
	Department of Combinatorics and Optimization\\
	University of Waterloo\\
	\tt pu.gao@uwaterloo.ca \and 
	Mikhail Isaev\thanks{Research supported by ARC  DE200101045, ARC DP190100977.}\\
	School of Mathematics\\
	Monash University\\
	\tt mikhail.isaev@monash.edu \and
	Brendan D. McKay\thanks{Research supported by ARC DP190100977}\\
	School of Computing\\
	Australian National University\\
	\tt brendan.mckay@anu.edu.au}

\date{}

\begin{document}
	
	\maketitle

	\begin{abstract}
		Kim and Vu made the following conjecture [\textit{Advances in Mathematics}, 2004]: if $d\gg \log n$, then the random $d$-regular graph $\G(n,d)$ can be ``sandwiched'' between $\G(n,p_*)$ and $\G(n,p^*)$ where $p_*$ and $p^*$ are both asymptotically equal to $d/n$.
		This famous conjecture was previously proved for all $d\gg (n\log n)^{3/4}$.  
		In this paper, we confirm the conjecture when  $d \gg \log^4 n$. We also extend  this result to near-regular degree sequences.
	\end{abstract}
	
	\section{Introduction}
	The binomial random graphs, introduced by Erd\H{o}s and R\'{e}nyi~\cite{ER,ER2} around 1960,  are the best known and most thoroughly studied among random graph models.
In the binomial random graph $\G(n,p)$, the vertex set is $[n]=\{1,2,\ldots,n\}$, and between each pair of vertices an edge is present independently with probability $p$.
More recently, the random regular graph $\G(n,d)$, and more generally $\G(n,\bfd)$, the uniformly random graph with degree sequence $\bfd=(d_1,\ldots,d_n)$, have received much attention.
In studies of random regular graphs, it was noticed that $\G(n,d)$ shares many properties with $\G(n,p)$ when $d\gg \log n$, where $p\approx d/n$.
For instance, asymptotically, they have the same chromatic number, the same diameter,  the same independence number, and the same order of the spectral gap.
This motivated the following sandwich conjecture, proposed by Kim and Vu~\cite{KimVu} in 2004.
	
	\begin{Conjecture}[Sandwich conjecture]\label{con:sandwich}
		For every $d\gg \log n$, there are $p_*=(1-o(1))d/n$, $p^*=(1+o(1))d/n$ and a coupling $(G_*, G, G^*)$ such that $G_*\sim \G(n,p_*)$, $G^*\sim \G(n,p^*)$, $G\sim \G(n,d)$ and $\Pr(G_*\subseteq G\subseteq G^*)=1-o(1)$.
	\end{Conjecture}
	
	The condition $d\gg \log n$ in the conjecture is necessary. When $d=O(\log n)$, there is a constant $\eps>0$ such that, with   probability at least $\eps $, there exist two vertices in $\G(n,p=d/n)$  which  degrees differ by at least  $\eps d$.
	Therefore, Conjecture~\ref{con:sandwich} cannot hold for this range of~$d$. For $\log n\ll d\ll n^{1/3}/\log^2 n$,
	Kim and Vu proved a weakened version of the sandwich conjecture where $G_* \subseteq G$ holds, but  $G\subseteq\up$ is replaced by a bound on $\varDelta(G\setminus \up)$ (see the precise statement in~\cite[Theorem 2]{KimVu}\footnote{Vu has confirmed that $\varDelta(\up\setminus G)$ in their theorem is a typo for $\varDelta(G\setminus \up$).}). 
	
	The partial sandwiching result of Kim and Vu \cite{KimVu} was further extended to $d=o(n)$ by Dudek, Frieze, Ruci{\'n}ski  and {\v{S}}ileikis~\cite{Dudek}, with basically the same coupling scheme used by Kim and Vu, and refinement. The first full sandwich result (i.e.\ $G_*\subseteq G\subseteq G^*$) was obtained by the authors  for $\min\{d,n-d\}\gg n/\sqrt{\log n}$ in~\cite{gao2020,gao2020Sandwiching}. Finally, Klimo\v{s}ov\'{a},   Reiher,  Ruci\'{n}ski and  \v{S}ileikis~\cite{Klimosova} confirmed that Conjecture~\ref{con:sandwich} holds for all  $\min\{d,n-d\}\gg (n\log n)^{3/4}$.
	
	One of the main contributions of this paper is confirmation of the sandwich conjecture
 for  $\min\{d, n- d\}\gg \log^4 n$. 
	\begin{theorem}\label{thm:sandwich}
		Conjecture~\ref{con:sandwich} holds for all
		$d$ such that $\min\{d,n-d\}\gg \log^4 n$.
	\end{theorem}

	\begin{Remark}\label{r:gap}
		Given Theorem~\ref{thm:sandwich} there are only two ranges where the Kim--Vu sandwich conjecture remains undecided:
		(i)~$ \log n\ll d = O(\log^4 n )$
		and (ii)~$d =n - O(\log^4 n)$.
		\begin{itemize}
			\item[(a)]  We have tried  to optimise choices of various parameters in the proof to obtain the best possible range of $d$, but  some new additional ideas seem to be necessary to close gap~(i).
		
			\item[(b)] Gap (ii) is simpler than gap~(i) because the complete graph (which is $G(n,p)$ with $p=1$) can be used as the top side of the sandwich.  Our coupling procedure with a few variations works to establish Conjecture~\ref{con:sandwich} in this range; see \cite[Theorem~1.5(c)]{gao2020arxiv} and its proof in our preprint~\cite{gao2020arxiv}. However, we do not include it in this paper to avoid unnecessary repetitions and to improve the clarity of the presentation.
		\end{itemize}
	\end{Remark}
	
	\smallskip
	
	The sandwich theorem not only explains the universality phenomenon between binomial random graphs and random regular graphs, it also gives a powerful tool for analysing random regular graphs.
Although $\G(n,d)$ has been extensively studied, much less has been proved in comparison to~$\G(n,p)$.
Proofs for properties of $\G(n,d)$ where $d=\omega(1)$ typically involve highly technical analysis, as it is not easy to estimate subgraph probabilities.
One would use the switching method to estimate such probabilities when $d=o(n)$, and for $d=\Theta(n)$, enumeration arguments involving multi-dimensional complex integrals are usually applied.
In the range of $d$ where the sandwich conjecture is confirmed, one can immediately translate monotone or convex properties, phase transitions, and many graph parameters from $\G(n,p)$ to $\G(n,d)$ with little effort.
See  examples of such translations in~\cite[Section 3]{gao2020arxiv}.

	A stronger version of Conjecture~\ref{con:sandwich} is given in~\cite[Conjecture 4]{gao2020Sandwiching} which extends Kim--Vu's sandwich conjecture to $\G(n,\dvec)$ when $\dvec$ is near-regular. We state the conjecture below.
	
	Given   ${\dvec}=(d_1,\ldots, d_n) \in \Reals^n$,
	let $\range{\dvec}$ stand for the difference
	between the maximum and minimum components of~$\dvec$. Denote 
	$\varDelta(\dvec) = \max_j d_j$.
	If  $\degvec(G)$ is the degree sequence of  a graph $G$, we will also use notations $\varDelta(G) = \varDelta(\degvec(G))$ and $\range{G} = \range{\degvec(G)}$.  
	
	\begin{Definition}  
		A  sequence $\dvec=\dvec(n) \in \{0,\ldots, n-1\}^n$  is called near $d$-regular if 
		\begin{align*}	
			|\varDelta(\dvec)-d| \ll  \min\{d,n-d\} \qquad \text{and} \qquad 
			\range{{\dvec}}  \ll  \min\{d,n-d\}
    		\end{align*}
		as $n\rightarrow \infty$, where $d=d(n)$ is a sequence of positive real numbers. 
	\end{Definition}

	\begin{Conjecture}\label{con:nonregular}
		Assume ${\dvec}= \dvec(n)$ is a near $d$-regular degree sequence
		with $d\gg \log n$. Then there are $p_*=(1-o(1))d/n$ and $p^*=(1+o(1))d/n$ and a coupling $(\low, G, \up)$ such that $\low\sim \G(n,p_*)$, $\up\sim \G(n,p^*)$, $G\sim \G(n,{\dvec})$ and $\Pr(\low\subseteq G\subseteq \up)=1-o(1)$.
	\end{Conjecture}

	\begin{Remark}
	We remark that the conditions cannot be weakened to  
 \[ 
 \range\dvec = o(\varDelta(\dvec)),
 \] 
 that is,  only requiring all degrees  to be asymptotically equivalent ignoring similar restrictions for the complementary degree sequence. Here we give an example illustrating this for the case 
  when $d = n -o(n)$.
Since we can  take $G^*$ to be the complete graph,  the existence of a coupling  $(\low, G, \up)$  is equivalent to embedding $\G(n,1-o(1))\subseteq \G(n,\dvec)$.
Let $\dvec$ be such that its complement degree sequence follows a power law with exponent between 2 and 3. It is implied by~\cite[Lemma 3]{Gao18} that for such $\dvec$, there exist pairs of vertices for which the edge probability between them is $o(1)$. Hence, it is not possible to embed $\G(n,1-o(1))$ into $\G(n,\dvec)$. More examples of co-sparse $\dvec$ with $\range{\dvec} = o(n)$ which do not allow such an embedding can be found in~\cite{Gao16}.
% We expect that the condition $\range\dvec\ll\min\{d,n-d\}$ in Conjecture~\ref{con:nonregular}
% can be weakened in the case $n-d=o(n)$.
	\end{Remark}

Conjecture \ref{con:nonregular} holds for $\min\{d,n-d\} = \Theta(n)$  by \cite[Theorem 5]{gao2020Sandwiching}.
	 The next result confirms  this conjecture  for all ``sufficiently'' near $d$-regular sequences.
	
	\begin{theorem}
		\label{thm:sandwich2}
		Conjecture~\ref{con:nonregular} is true for all  $\dvec =\dvec(n)$  satisfying
		\[ 
		\min\{d, n-d\}\gg \log^4 n \qquad \text{and} \qquad
		\range{\dvec} \ll   \dfrac{1}{\log n}\min\{d, n-d\}.\]
	\end{theorem}
	Note that Theorem~\ref{thm:sandwich} now follows as an immediate corollary of Theorem~\ref{thm:sandwich2}.
	We prove Theorem \ref{thm:sandwich2} in Section \ref{s:proof_sandwich}.

	\smallskip
	
	Theorem~\ref{thm:sandwich2} does not solve Conjecture~\ref{con:nonregular} completely for the reasons explained in Remark \ref{r:gap}. 
	Also, note that we  assumed a tighter bound on $\range{\dvec}$ in Theorem~\ref{thm:sandwich2} than in Conjecture~\ref{con:nonregular}.  Thanks to~\cite[Theorem 5]{gao2020Sandwiching}, we can restrict our attention to $d=o(n)$, since the case $n-d=o(n)$ follows by taking complements.
	 The following result establishes  the bottom side of the sandwich   for the case of $d= o(n)$.

	\begin{theorem}\label{thm:low}
		Assume ${\dvec}= \dvec(n)$ is a near $d$-regular degree sequence
		with  $\log n\ll d\ll n$. Then there are $p_*=(1-o(1))d/n$ and a coupling $(\low, G)$ such that $\low\sim \G(n,p_*)$, $G\sim \G(n,{\dvec})$ and $\Pr(\low\subseteq G) = 1-o(1)$.
	\end{theorem}
	Theorem \ref{thm:low} is a simplification of Theorem~\ref{thm:bottom}  stated and proved in Section \ref{s:bottom}.

	Throughout the paper, we assume that $\dvec$ is a realisable degree sequence, i.e.\ $\G(n,\dvec)$ is nonempty. This necessarily requires that  ${\dvec}$ has nonnegative integer coordinates and even sum.  All asymptotics in the paper refer to $n\to\infty$. For two sequences of real numbers $a_n$ and $b_n$, we say $a_n=o(b_n)$ if $b_n\neq 0$ eventually and $\lim_{n\to\infty} a_n/b_n=0$. We say $a_n=O(b_n)$ if there exists a constant $C>0$ such that $|a_n|\leq C\,|b_n|$ for all $n$.
	We write $a_n=\omega(b_n)$ or $a_n=\Omega(b_n)$  if $a_n>0$ always
	and $b_n=o(a_n)$ or $b_n=O(a_n)$, respectively.
	If both $a_n$ and $b_n$ are positive sequences, we will also write
	$a_n\ll b_n$ if $a_n=o(b_n)$, and $a_n\gg b_n$  if $a_n=\omega(b_n)$.

    For a vector $\xvec \in \Reals^n$ and $\rho\geq 1$ we use the standard norm notation
    \[
    \norm\xvec_\rho := \biggl(\sum_{i=1}^n \abs{x_i}^\rho\biggr)^{\!1/\rho}
    \text{~~and~~}
    \norm \xvec_\infty := \max_{i=1}^n \nolimits\,\abs{x_i}
\]
For an $n\times n $ matrix $M$ and $\rho \in [1,\infty]$, 
we also use the induced matrix norm  defined by
\[
     \norm M_\rho := \max_{\xvec \neq \boldsymbol{0} }\dfrac{\norm{M\xvec}_\rho}{\norm\xvec_\rho}.
\]
All of these norms are subadditive ($\norm{\xvec_1+\xvec_2}_\rho
\leq \norm{\xvec_1}_\rho+\norm{\xvec_2}_\rho$ and
$\norm{M_1+M_2}_\rho
\leq \norm{M_1}_\rho+\norm{M_2}_\rho$)
and submultiplicative ($\norm{M_1 M_2}_\rho
\leq \norm{M_1}_\rho\norm{M_2}_\rho$ and
$\norm{M\xvec}_\rho
\leq \norm M_\rho\norm\xvec_\rho$).

\subsection{Proof ideas}
    
	The coupling procedures used in all previous work~\cite{KimVu,Dudek,gao2020Sandwiching,Klimosova} are essentially the same (for simplicity we explain the coupling procedure for $\G(n,d)$ only, and the treatment for $\G(n,\dvec)$ is similar): a triple of random graphs $(F_*,F,F^*)$ are sequentially constructed simultaneously so that every time an edge is added to $F$ according to the correct marginal distribution (of $\G(n,d)$ containing that edge conditional on the current construction), it is added to the bottom graph $F_*$ with a slightly smaller probability, independent of the previous construction.
When~$F$ is close to fully constructed, the construction of $F$ is completed without changing $F_*$ and~$F^*$. The coupling output is $(F_*,F,F^*)$.
By choosing several parameters carefully, $F_*$ and $F^*$ are binomial random graphs with asymptotically equal edge density, and $F_*\subseteq F$ with probability $1-o(1)$, but $F$ is not a subgraph of~$F^*$.
This idea was introduced by \cite{KimVu} and then extended to the hypergraph setting in \cite{Dudek}. 
For more details about this coupling; see Section \ref{s:bottom}.
	
	To overcome the problem that  $F$ is not a subgraph of $F^*$, the authors~\cite{gao2020Sandwiching} used a similar coupling construction for embedding $\G(n,1-p_*)$ inside $\G(n,n-d-1)$ in the dense setting. Complementing this coupling gives the top side of the sandwich. This idea was also used in the sandwich results of \cite{Klimosova} to allow degrees where $\min\{d,n-d\}\gg (n\log n)^{3/4}$.

	In this paper, we propose a novel 2-round coupling scheme. First,  we construct $(F_*,F,F^*)$ as before but stop when $F$ is sufficiently close to completion. Then, in the second round of coupling, we find a random subgraph $F'$ of $K_n\setminus F$ with the correct distribution, so that $F\cup F'$ has the distribution $\G(n,d)$. We do so by deleting  edges from the complement of $F$ according to the right marginal distribution until reaching a $d$-regular graph $G$. Since at every step edges are added or deleted from the complement according to the correct marginal distribution, we have that $G\sim \G(n,d)$.  The edge deletion process for $G$ is coupled with another edge deletion process $G'$ such that
	$G'$ and $F^*$ are two independent binomial random graphs and $G\subseteq G'\cup F^*$
	with probability $1-o(1)$. This yields a coupling that embeds $G\sim \G(n,d)$ inside $G'\cup F^*\sim\G(n,p)$ for some~$p$. We describe the second round of the coupling scheme in detail in Section \ref{s:top1}.
	
	A spanning subgraph with degree sequence $\tvec = (t_1,\ldots,t_n)$ of a graph $S$ is called $\tvec$-factor of $S$. Let $S_{\tvec}$ denote a uniformly random $\tvec$-factor of $S$. 
	The analysis of our new  2-round coupling procedure requires   that all edges  are approximately  equally likely to appear in $S_{\tvec}$ for a near-regular sequence $\tvec$ and a pseudorandom graph $S$.
	Such results were previously known only for dense $S$  but not for sparse $S$; see~\cite[Theorem 10]{gao2020Sandwiching} and references therein.
	To allow sparse $S$, we use a switching method that involves counting walks in $S$ up to logarithmic length, which to our knowledge has not been done before. Thus, the following theorem is of independent interest.

 Let $J$ denote the $n\times n$ matrix with all entries equal to $1$.  Given a graph $S$, let 
$A(S)$ denote its adjacency matrix and let $|S|$ be its number of edges. Throughout the paper, we let $N$ denote $\binom{n}{2}$, which is the number of edges in a complete graph on $n$ vertices.

% \mi{How about this one? Take $\alpha = -2 \log c/ \log n$}

\begin{theorem}\label{l:co-sparse:simple}
Let $c \in (0,1)$ be fixed and  $S = S(n)$ be  graphs on vertex set $[n]$ such that  
\begin{itemize}
\item[\rm (i)] $\range{S} = o(\varDelta(S)/\log n)$; 
\item[\rm (ii)]
$ \|A(S) - \frac{|S|}{N}J\|_2 \leq c \varDelta(S) $;
\item[\rm (iii)]
for any fixed $\delta>0$, we have
			\[
\sum_{(jk)\st jk \in S} x_j y_k = (1+o(1)) \|\xvec\|_1\|\yvec\|_1 \dfrac{\varDelta(S)}{n}.
			\]
            uniformly over  all $\xvec,\yvec \in [0,1]^n$  with $\|\xvec\|_1, \|\yvec\|_1  \geq \delta n$.
\end{itemize}
Let   $\tvec=\tvec(n) \in \Naturals^n$  be degree sequences satisfying 
\begin{itemize}
\item[\rm (iv)] $\range{\tvec} = o(\varDelta(\tvec)/\log n)$ and $\dfrac{\log^2 n}{\varDelta(S)} \ll \varDelta(\tvec)  \ll \dfrac{\varDelta(S)}{\log n}$.
\end{itemize}
 Then,  uniformly over $jk\in S$, 
 \[
        \Pr(jk \in S_{\tvec}) = (1+o(1)) \frac{\varDelta(\tvec)}{\varDelta(S)}.
 \]
\end{theorem}

Theorem~\ref{l:co-sparse:simple} is a corollary of Theorem~\ref{l:co-sparse} by taking $\alpha= -\frac12\log c/\log n$ and $\beta = 2 $. The proof of the latter theorem will be the topic of Section~\ref{s:switchings}. 

The restriction of Theorem \ref{l:co-sparse:simple} on the relative density $\varDelta(\tvec)/\varDelta(S) \ll 
\dfrac{1}{\log n}$ is the fundamental reason determining our range $d\gg \log^4 n$.
Any improvement on this would allow sparser degree sequences, in particular,  an analog of Theorem \ref{l:co-sparse:simple} under a relaxed condition $\varDelta(\tvec)/\varDelta(S) = o(1)$ would resolve Conjecture \ref{con:sandwich} completely.

	\section{Sandwich}
	This section is organised as follows. We treat the bottom and top sides of the sandwich separately in Section \ref{s:bottom} and Sections \ref{s:top1}, \ref{s:top2}, respectively. Then we prove Theorem \ref{thm:sandwich2}. 
	For both the bottom and top sides,
	our coupling construction heavily relies on  edge probability estimates for a random near-regular $\tvec$-factor in a pseudorandom graph $S$. 
	The bottom side is more straightforward as it only requires enumerating sparse factors in an almost complete graph, which was done by McKay \cite{McKay1981}.
	For the top side, we use our new estimate stated in the previous section as Theorem \ref{l:co-sparse:simple}.

	\subsection{Bottom of the sandwich}\label{s:bottom}
	
	In this section, we construct a coupling $\G(n,p) \subseteq \G(n,\dvec)$.  The next lemma is crucial for our construction.

	\begin{lemma}\label{l:sparse}
		Let  $S$ be a graph on  $n$ vertices and   $\tvec = (t_1,\ldots,t_n)$ be a degree sequence such that
		the set of $\tvec$-factors of $S$ is not empty and 
		\[
		\varDelta (\tvec) (\varDelta(\tvec) +\varDelta(K_n\setminus S))  \ll tn, \qquad \text{where }\  t= \frac{t_1+\cdots +t_n}{n}>0.
		\]
		Then, for all $jk \in S$, we have
		\[
		\Pr(jk\in S_{\tvec}) = \left(1+ O\left(\frac{\varDelta (\tvec) (\varDelta(\tvec) +\varDelta(K_n\setminus S))}{tn}\right)\right)\frac{t_j t_k}{t n}.
		\] 
	\end{lemma}
	\begin{proof}
		We follow the notation in~\cite{McKay1981} and apply the bounds of \cite[Corollary 2.4]{McKay1981} and  \cite[Lemma 2.8]{McKay1981} with  $\gvec = \tvec$, $H =\emptyset$, 
		and $L = (K_n \setminus S) \cup \{jk\}$ to estimate the ratio $N(\gvec, L, jk)/N(\gvec, L,\emptyset)$,
		where $N(\gvec,L,A)$ denotes the number of graphs $H$ with degree sequence $\gvec$ such that $H\cap L=A$.
		Observing that  
		\[
		\Pr(jk\in S_{\tvec})
		=  \frac{N(\gvec, L, jk)}{N(\gvec, L, jk)+N(\gvec, L,\emptyset)},
		\]
	  the proof is complete.
	\end{proof}

	Let $G_{\dvec}$ be a random graph distributed according to  
	$
	\G(n,\dvec) 
	$
	and 
	\begin{equation*}
M:=|G_{\dvec}|= \dfrac12 \sum_{i\in [n]}d_i.% \label{eq:N-M}
	\end{equation*} Recall that $N:= \binom{n}{2}$.
	Given $0\leq m\leq M$, let
	$\G (n,\dvec, m)$ be the model of random graph obtained from 
	$G_{\dvec} $ by choosing $m$ of its edges uniformly at random.  
	For a random graph $H$, we say $H \sim \G(n, |H|)$ to mean that its conditional distribution on the event $|H|=m$ (for any eligible $m$) is uniform.
      Let $\Bin(a,q)$ denote the binomial distribution with parameters $a,q$ where $a$ is a positive integer and $q\in [0,1]$.
    % %or $F \sim \G(n,\dvec,|F|)$ 
    % to mean that
    % $F$ is distributed as $\G(n,X)$ 
    % %or $\G(n,\dvec,X)$, respectively,
    % where $X$ is the random variable $|F|$.
	
	\begin{theorem}\label{thm:bottom}
		Let $\dvec=\dvec(n)$ be a  near $d$-regular  sequence, where $d=d(n)=o(n)$, 
		and 
		$\xi = \xi(n) = o(1)$ are positive sequences such that  
		\begin{equation}
			\range{\dvec} \leq \xi d \qquad
			\text{and}
			\qquad
			\xi n \geq d \gg \xi^{-3}\log n. \label{eq:d-range}
		\end{equation}
		Then there exists a coupling $(G_*, G)$ with 
		$G_* \sim \G(n, p_*)$ and   
		$G \sim \G(n,\dvec)$  
		such that  
		\[
		p_* = (1+O(\xi))\dfrac{d}{n}, \qquad \Pr(G_* \not\subseteq G) \leq e^{-\omega(\log n)}. 
		\]
	\end{theorem}
	
	\begin{proof}
		Let $\zeta=C\xi$ where $C>0$ is a sufficiently large constant (which depends only on the implicit constant in the $O(\,)$ bound of Lemma \ref{l:sparse}).
        	Let 
		\begin{equation*} 
			\Time \sim \po(\mu), \qquad  \mu:=(1-\xi)M.
		\end{equation*}  
		We will show that the pair $(G_*,G):= (G_{\zeta}^{(\Time)}, G^{(\infty)})$ defined below satisfies the theorem requirements.
		
		Consider a sequence $(a^{(\time)})_{\time \in \Naturals}$ of   random numbers
		and a sequence $(e^{(\time)})_{\time \in \Naturals}$ of random edges 
		chosen uniformly and independently from $[0,1]$  and $K_n$, respectively.
		Define $G_\zeta^{(0)} = G^{(0)} = G_0^{(0)} := \emptyset$. For $\time \geq 1$, let
		\begin{align*}
			G_0^{(\time)}&:= G_0^{(\time -1)} \cup e^{(\time)};\\
			G_{\zeta}^{(\time)} &:= 
			\begin{cases}
				G_\zeta^{(\time-1)} & \text{if $a^{(\time)}\ge1-\zeta$,}\\[0.5ex]
				G_\zeta^{(\time -1)} \cup e^{(\time)} &\text{otherwise;}
			\end{cases}
			\\[0.6ex]
			G^{(\time)} &:=
			\begin{cases}
				G^{(\time-1)} & \text{if $a^{(\time)}\ge1-\eta^{(\time)}
					$,
				}\\
				G^{(\time -1)} \cup e^{(\time)} &\text{otherwise,}
			\end{cases}
		\end{align*}
		where  $\eta^{(\time)}=\eta(G^{(\time-1)},e^{(\time)})$
		\begin{equation}\label{def:eta}
			\eta(G^{(\time-1)},e^{(\time)}):= 
			\begin{cases}
				1 &\hspace{-9em}
				\text{if   $e^{(\time)} \in G^{(\time-1)}$ or $|G^{(\time-1)}|=M$,}
				\\[0.5ex]	 
				\displaystyle
				1 - \frac{\Pr\left(e^{(\time)} \in G_{\dvec}\mid 
					G^{(\time-1)} \subseteq G_{\dvec}\right)}
				{\max_{e \notin G^{(\time-1)} }\Pr\left(e \in G_{\dvec}\mid 
					G^{(\time-1)} \subseteq G_{\dvec}\right)} & \text{otherwise.}
			\end{cases}
		\end{equation}

		Note that if $G^{(\time-1)}$ has fewer than $M$ edges
		then 
		$\Pr\left(  G^{(\time-1)} \subseteq G_{\dvec}\right)>0$,
		and so there is an edge $e \notin G^{(\time-1)} $
		such that $\Pr\left(e \in G_{\dvec}\mid 
		G^{(\time-1)} \subseteq G_{\dvec}\right)>0$.
		Therefore, $G^{(\time)}$ eventually becomes a graph with degree
		sequence $\dvec$, which we denote by $G^{(\infty)}$.

		Note that the graphs $G_0^{(\time)}$ and $G_\zeta^{(\time)}$
		are identical to the underlying graphs of the multigraphs 
		$M_0^{(\time)}$ and $M_\zeta^{(\time)}$ in the procedure $Coupling()$ from \cite[Figure 1]{gao2020Sandwiching}.
        Our definition of $G^{(\time)}$
		is slightly different from  \cite[Figure 1]{gao2020Sandwiching}---here
		we avoid calling the procedure $IndSample()$ so that $G^{(\time)}$ is  well-defined for every $\time$.
		First, we list all the bounds from \cite{gao2020Sandwiching} that we use in this proof. 
		Let 
		\begin{equation*}%\label{def:m-whtout-hat}
		m^{(\time)} = |G^{(\time)}|,\qquad 	m_0^{(\time)} = |G_0^{(\time)}|,
		\qquad 	m_\zeta^{(\time)} = |G_\zeta^{(\time)}|. 
		\end{equation*}
Also, by  \cite[Lemma 2, Lemma 7]{gao2020Sandwiching},
		\[
		G_0^{(\time)}\sim \G(n,m_0^{(\time)}), \qquad  G_\zeta^{(\time)}\sim \G(n,m_\zeta^{(\time)})
		\]
and, for the Poisson random step $\time = \Time$,  we have
		\[
		G_0^{(\Time)}  \sim \G(n,p_0), \qquad G_\zeta^{(\Time)}  \sim \G(n,p_\zeta),
		\] 
		where
		\[
		p_0:=1 - e^{-\mu/N}, \qquad p_\zeta:=1 - e^{-(1-\zeta)\mu/N}.
		\]
		By the  assumptions, we find that  $M/N = (1+O(\xi))d/n$  and thus 
        \[ % begin{equation}\label{eq:p-value}
         p_\zeta,p_0 = (1+O(\xi)) d/n   
        \]
        % \jc{was --- $p_\zeta = (1+O(\xi)) d/n$} as required \jt{(recalling that $\zeta$ is set $C\xi$ for some constant $C>0$).}  
        Clearly, by our construction,  we have $G_\zeta^{(\Time)} \subseteq G_0^{(\Time)}$.

During the refereeing process of this paper, it was pointed out by an anonymous referee that \cite[Lemma 3]{gao2020Sandwiching} requires more careful consideration.
It claims that 
		\begin{equation}\label{distr_G}
			G^{(\time)} \sim \G (n,\dvec, m^{(\time)}),
		\end{equation}  
which is wrong. Fortunately, for our purposes it is sufficient that the claimed distribution \eqref{distr_G} is correct   at steps $\tau$ where we increase the number of edges in $G^{(\time)} $, and this is exactly what we get from the arguments in the proof of \cite[Lemma 3]{gao2020Sandwiching}. 
We state it below for ease of reference later and give a detailed proof in Section \ref{S:crucial}. The results and the proofs in~\cite{gao2020Sandwiching} are unaffected after replacing~\cite[Lemma 3]{gao2020Sandwiching} by the following claim.
\begin{claim}\label{claim:crucial}
If  $i\in [M]$ and $\time_i$ is the first step where $|G^{(\time)}|$ has $i$ edges  then $G^{(\time_i)}\sim \G (n,\dvec, i).$ 
\end{claim}
\begin{Remark}
One might think that, since $G^{(\tau)}$ does not change between the steps $\tau_{i}$ and $\tau_{i+1}$, it should maintain the  distribution as claimed in \eqref{distr_G}.   A subtle issue is that the ``waiting time" $\tau_{i+1}-\tau_{i}$, in fact depends on the structure of $G^{(\tau_i)}$ and  variations between the conditional edge probabilities in \eqref{def:eta}.
\end{Remark}

		In particular, it follows from Claim \ref{claim:crucial}  that $G^{(\infty)} \sim \G(n,\dvec)$.  
		Thus, it remains to prove  the claimed probability bound for the event 
		$G_{\zeta}^{(\Time)} \not\subseteq G^{(\infty)}$. Note that, for our choice of $\mu$,
		\[
		p_0=  1- e^{-(1-\xi)M/N} \leq (1-\xi)M/N.
		\]
		Applying the Chernoff bound for $ |\G(n,  (1- \xi)M/{N})|\sim\Bin(N,(1-\xi)M/N)$ (see  Lemma~\ref{l:technical}),  we get that  
		\begin{equation}\label{eq:lowM}
			M- m^{(\Time)} \geq M-m_0^{(\Time)} 
			\geq M - |\G(n,  (1- \xi)M/{N})|
			\geq \xi M/2
		\end{equation}
		with probability at least $1- e^{-\Omega(\xi M)} \geq 1 -e^{-\omega(\log n)}$.
		% Using Lemma \ref{l:technical} again for $\Time\sim\po(\mu)$, we estimate that 
		% \[
		% \Pr(\Time > 2\mu ) = e^{-\Omega(\mu)}.
		% \]
By our construction, we have  $G^{(\Time)}\subseteq G^{(\infty)}$, therefore
$G_{\zeta}^{(\Time)} \not \subseteq G^{(\infty)} $ implies that 
$G_{\zeta}^{(\Time)} \not \subseteq G^{(\Time)} $. Moreover,  if $G_{\zeta}^{(\Time)} \not \subseteq G^{(\Time)} $ occurs then one of the two events must occur:
\begin{itemize}
    \item $m^{(\Time)} > (1-\xi/2) M$,
    \item $\eta(G^{(\tau_i)},e)>\zeta$ for some $e\notin G^{(\tau_i)}$ and   $i\leq (1-\xi/2) M$.  
\end{itemize}
Note that we only need to consider steps $\tau_i$ since  $G^{(\tau)} = G^{(\tau_i)}$ for all $\tau_i\leq \tau<\tau_{i+1}$.
By the union bound,
     \begin{equation}\label{eq:union:sparse}
     \begin{aligned}
			\Pr(G_{\zeta}^{(\Time)} \not \subseteq G^{(\infty)} ) &\leq
			\Pr(G_{\zeta}^{(\Time)} \not\subseteq G^{(\Time)} )
			\\
            &\leq   e^{-\omega(\log n)}   
            +   \sum_{i \leq (1-\xi/2)M} \Pr\left( \max_{e\notin G^{(\tau_i)} }\left\{\eta(G^{(\tau_i)},e)\right\} >\zeta\right).
        %    \\
            % &
            % =e^{-\Omega(\mu)} +  2\mu  \max_{i \in [M]\st \time_i+1 \leq 2\mu} \, \Pr\left(\bigl(\eta^{(\time_i+1)} >\zeta\bigr)\cap \bigl(e^{(\time_i+1)}\notin G^{(\time_i)}\bigr)\right)
            \end{aligned}
		\end{equation}

	Next, for each $i \leq (1-\xi/2)M$, we estimate the probability that there is $e\notin G^{(\tau_i)}$  such that $\eta(G^{(\tau_i)},e)>\zeta$. 
%$\Pr\left(\bigl(\eta^{(\time)} >\zeta\bigr)\cap \bigl(e^{(\time)}\notin G^{(\time-1)}\bigr)\right)$ 
By Lemma~\ref{l:sparse} with
		\[
		S:= K_n \setminus G^{(\time_i)}, \qquad
		\tvec := \dvec - \gvec^{(\time_i)},
		\] where
%\[
 %   \time' = \max_{i \in [M]}\{\time_i \st \time_i < \time\}
%\] 
%and 
$\gvec^{(\time_i)}$ is the degree sequence of $G^{(\time_i)}$. 
        By~Claim \ref{claim:crucial},   we have $G^{(\tau_i)}\sim \G(n,\dvec,i)$. Thus we get that, for all $j\in [n]$,
		\begin{equation}\label{t-hyper}
			t_j \sim \operatorname{Hypergeometric}(M, d_j, M-i).
		\end{equation}
		Indeed, for every graph $X\in \G(n,\bfd)$, after sampling a uniformly random set of $M-i$ edges from the set of total $M$ edges of $X$, $t_j$ counts the number of edges incident with $j$ that are sampled, out of the total $d_j$ edges incident with $j$ in $X$.
        Note that
        \[
            \E \tvec = \dfrac{M-i}{M} \dvec. 
        \]
		Applying the Chernoff bound for all $t_j$ (see Lemma \ref{l:technical} with $\eps = \xi$) and using \eqref{eq:lowM}, we get that  
		\begin{align*} 
			\Pr \left(	\left\| \tvec - \E \tvec \right\|_\infty  >
			\xi \dfrac{M-i}{M}\,   d  \right)
            	 \leq  n e^{- \Omega(\xi^3 d)}
			= e^{- \omega(\log n)}.
		\end{align*}
		% For the second last inequality above we used
		% \[
		% \E [\tvec \mid m^{(\time_i)}=m] = \dfrac{M-m}{M} \dvec
		% \] 
		% and
		Since $\dvec$ is near $d$-regular, we have $M = (1+o(1))nd/2$. Note also that 
        \[
		t = \dfrac{t_1+\cdots +t_n}{n} =  2\,\dfrac{M-i}{n}  \geq \xi M/n.
		\]
		Since $\range{\dvec} \leq \xi d$, we get that with probability at least $1-e^{- \omega(\log n)}$, 
		\begin{equation}\label{eq:range_t}
			\range{\tvec} \leq \range{ \E \tvec} + \| \tvec - \E \tvec  \|_\infty \leq \dfrac{M-i}{M}\range{\dvec} +  \xi \dfrac{M-i}{M}    d  =  O(\xi t).
		\end{equation}
		We  also bound 
		\[
		\varDelta(K_n \setminus S) = \varDelta(G^{(\time_i)}) = O(d).
		\]
		Applying Lemma \ref{l:sparse} and taking  $C$ in the definition of $\zeta$ to be sufficiently large,  we find that
		\begin{equation}\label{eq:finally}
			\max_{e\notin G^{(\tau_i)} }\left\{\eta(G^{(\tau_i)},e)\right\} =O\left( \dfrac{\range{\tvec}}{\varDelta(\tvec) - \range{\tvec}} + \dfrac{\varDelta(\tvec)^2}{tn} + \dfrac{\varDelta(\tvec)\varDelta(K_n-S)}{tn}\right) 
			=   O(\xi) <\zeta
		\end{equation}
		with probability  at least $1-e^{- \omega(\log n)}$.

Finally, combining \eqref{eq:union:sparse} and  \eqref{eq:finally}, we get  
		\[
		\Pr(G_{\zeta}^{(\Time)} \not \subseteq G^{(\infty)} ) 
		\leq e^{-\omega(\log n)} + (1-\xi/2)M e^{-\omega(\log n)}
		=e^{-\omega(\log n)}
		\]
		as claimed.
	\end{proof}

    With a slight tweak of the arguments of  Theorem \ref{thm:bottom} we also establish  the following lemma, which will be useful for the top side of the sandwich. 
	
	\begin{lemma}\label{L:top-partial}
		Let the assumptions of Theorem \ref{thm:bottom} hold. % Let $\hat{m} = \lceil(1-5\xi)M\rceil$.
        Then there is 
        $m_* =  (1-O(\xi))M$ and  
        a triple of random graphs
		$(G_*, F_*, G_0)$ such that 
		\begin{enumerate}
			\item[(i)] $G_* \sim \calG(n,p_*)$, $G_0 \sim \calG(n,p_0)$,  
			$F_* \sim \calG(n,\dvec,m_*)$, and
			\[
			p_*, p_0 = (1+O(\xi)) \dfrac{d}{n};
			\] 
			\item[(ii)] with probability at least $1-e^{-\omega(\log n)}$, we have 
			\[
			G_* \subseteq F_* \subseteq G_0, \quad \varDelta(\dvec-\degvec(F))=\Theta(\xi d), \quad \range{\dvec-\degvec(F)}=O(\xi^2 d).
			\]
			%where $\fvec$ is the degree sequence of $F_*$. 
		\end{enumerate}
	\end{lemma}
	\begin{proof}
	We use the same construction for thhe sequence of graphs $(G_0^{(\tau)}, G_\zeta^{(\tau)}, G^{(\tau)})$  as in the proof of  Theorem \ref{thm:bottom}.
   As before, let $\mu = (1-\xi)M$.  Recall also that 
    \[
         p_\zeta =  1 - e^{-(1-\zeta)\mu/N} =  (1+O(\xi))d/n.
    \]
    We take $\mu_* = (1-O(\xi))M$ and $m_* =(1-O(\xi))M$ to be any such that 
    \begin{equation}\label{def:hatm-mustar}
         (1+\xi )(1 - e^{-\mu_*/N})  \leq m_*/N \leq (1-\xi) p_{\zeta}.
    \end{equation}
        Let 
    \[
        \Time = \Time_* + \Time_{small}, 
    \]
    where $ \Time_*$ and $\Time_{small}$ are independent and 
    \[
        \Time_* \sim  \po(\mu_*), \qquad   \qquad  \Time_{small} \sim  \po(\mu-\mu_*). 
    \]
In particular, we have $\Time\sim\po(\mu)$ as before. Let 
$\Time'$ denote the first step $\tau$ when $G^{(\tau)}$ contains  $m_*$ edges.
Then, we  define the required triple $(G_*,F_*,G_0)$   by 
\[
    G_* :=G_\zeta^{(\Time_*)}, \qquad F_*:= G^{(\Time')}, \qquad G_0 := G_0^{(\Time)}.
\]
 It remains to show that it satisfies the stated properties (i) and (ii).

As in the proof of Theorem \ref{thm:bottom}, we have
\[
   G_0 \sim \calG(n,p_0), \text{ where } p_0 = 1 - e^{-\mu/N} = (1+O(\xi))d/n. 
\]
 Since $\Time_*$ is also Poisson, we get 
 \[ 
      G_* \sim \calG(n,p_*), \text{ where } 
      p_* =  1 - e^{-(1-\zeta)\mu_*/N} = (1+O(\xi))d/n.
    \]
Using   Claim \ref{claim:crucial}, we have established the desired distribution of $F$, completing the proof of~(i).
Also, from the proof of  Theorem \ref{thm:bottom},  
with probability $1 - e^{-\omega(\log n)}$, we have
\begin{equation}\label{sandwich-middle}
    G_\zeta^{(\tau)} \subseteq G^{(\tau)} \subseteq G_0^{(\tau)}
\end{equation}
for all $\tau \leq \Time$.   
Recalling that  
 $
    G_\zeta^{(\Time)}  \sim \G(n,p_\zeta)$ and applying the Chernoff bound (see Lemma~\ref{l:technical}) together with \eqref{def:hatm-mustar}, we find that 
    \[
       \Pr\left(|G_\zeta^{(\Time)}|<m_*\right)  \leq e^{-\Omega(\xi^2 Np_\zeta)} =e^{-\omega(\log n)}.
    \]
Similarly, since $\Time_*$ is also Poisson, we have 
$ 
    G_0^{(\Time_*)}  \sim \G(n,1 - e^{-\mu_*/N})$. 
Applying the Chernoff bound (see Lemma~\ref{l:technical}) together with \eqref{def:hatm-mustar}, we find that 
    \[
       \Pr\left(|G_0^{(\Time_*)}|>m_*\right)  \leq \exp \Big(-\Omega\left(\xi^2 N(1 - e^{-\mu_*/N})\right) \Big) =e^{-\omega(\log n)}.
    \]   
If  events $|G_\zeta^{(\Time)}|\geq m_*$,
$|G_0^{(\Time_*)}|\leq m_*$, and  \eqref{sandwich-middle} occur then   $\Time' \in [\Time_*,\Time]$ by definition of $\Time'$ and we  get  
\[
    G_\zeta^{(\Time_*)} \subseteq G^{(\Time_*)}
    \subseteq G^{(\Time')} \subseteq G^{(\Time)} \subseteq G_0^{(\Time)}.
\]
This establishes 
\[
   \Pr( G_* \subseteq F_* \subseteq G_0) \geq 1 - e^{-\omega(\log n)}.
\]
The stated bounds for $  \varDelta(\dvec-\degvec(F))$ and   $\range{\dvec-\degvec(F)}$ follow by applying the Chernoff bound to the hypergeometric random variables in exactly the same way as in the proof of Theorem~\ref{thm:bottom}; see \eqref{t-hyper} and \eqref{eq:range_t}.
\end{proof}

	% % % % % % % % % % %
	% % % % % % % % % % %
	% % % % % % % % % % %
	% % % % % % % % % % %
	% % % % % % % % % % %
	% % % % % % % % % % %
	% % % % % % % % % % %

	\subsection{Top of the sandwich}\label{s:top1}

	In this section, we construct a coupling $ \G(n,\dvec) \subseteq \G(n,p)$ based on the claim that all edges typically appear in our sequential construction with asymptotically the same probability. 
	Recall that 
	$G_{\dvec}$ is  a random graph distributed according to  
	$
	\G(n,\dvec) 
	$
	and 
	\[ 
	N:= \binom{n}{2}, \qquad M:=|G_{\dvec}| = \dfrac12 \sum_{i\in [n]}d_i. \]
	Recall also that $\G (n,\dvec, m)$ is the model of  a random graph obtained from 
	$G_{\dvec} \sim \G(n,\dvec)$ by choosing $m$ edges uniformly at random.

Recall that Theorem~\ref{thm:bottom} requires $\dvec$ to be a positive near $d$-regular degree sequence  satisfying~\eqref{eq:d-range} for some $d=o(n)$ and $\xi=o(1)$. The theorem below shows that with degree sequences satisfying similar but slightly stronger conditions, the top part of the sandwich holds.

    \begin{theorem}\label{T:top}   Let $\dvec=\dvec(n)$ be a  near $d$-regular  sequence, where $d=d(n)=o(n)$, 
		and 
		$\xi = \xi(n) \ll 1/\log n$ are positive sequences such that  
		\[
			\range{\dvec} \leq \xi d \qquad
			\text{and}
			\qquad
			\xi n \geq d \gg \xi^{-3}\log n.
		\]
		Then there exists a coupling $(G,G^*)$  with  $G^*\sim G(n,p^*)$,
		$G \sim \G(n,\dvec)$ such that 
		\[
		p^* = (1+ o(1)) \dfrac{d}{n}, \qquad
		\Pr(G \not\subseteq G^*) \leq e^{- \omega(\log n)}.
		\]              
	\end{theorem}

\begin{proof}

	 Let $F$ be an arbitrary graph on $[n]$ with $\degvec(F) \leq \dvec$ and $\zeta \in [0,1]$ which we specify more precisely later. Consider a sequence $(a^{(\time)})_{\time \in \Naturals}$ of   random numbers
		and a sequence $(e^{(\time)})_{\time \in \Naturals}$ of random edges 
		chosen uniformly and independently from $[0,1]$  and $K_n$, respectively.
		Let $H_\zeta^{(0)} = H^{(0)} = H_0^{(0)} = K_n$. For $\time \geq 1$, consider the following graph sequences:
		\begin{align*}
			H_0^{(\time)}&:= H_0^{(\time -1)} \setminus e^{(\time)};\\
			H_{\zeta}^{(\time)} &:= 
			\begin{cases}
				H_\zeta^{(\time-1)} & \text{if $a^{(\time)}\ge1-\zeta$,}\\[0.5ex]
				H_\zeta^{(\time -1)} \setminus e^{(\time)} &\text{otherwise;}
			\end{cases}
			\\
			H^{(\time)} &:=
			\begin{cases}
				H^{(\time-1)} & \text{if $a^{(\time)}\ge1-\eta_F^{(\time)}$}\\
				H^{(\time -1)} \setminus e^{(\time)} &\text{otherwise,}
			\end{cases}
		\end{align*}
		where  $\eta_F^{(\time)}=\eta_F(H^{(\time-1)},e^{(\time)})$
		\[  % begin{equation}\label{def:etaF}
			\eta_F (H^{(\time-1)},e^{(\time)}) := 
			\begin{cases}
				1 &\hspace{-14.5em} 
				\text{if    $|H^{(\time-1)}|=M$ or $e^{(\time)} \in F$, or $e^{(\time)}\notin H^{(\time-1)}$}
				\\[1ex]	 
				\displaystyle
				1 - \frac{\Pr\left(e^{(\time)} \notin G_{\dvec}\mid 
					F\subseteq G_{\dvec}\subseteq H^{(\time-1)}  \right)}
				{\max_{e \in H^{(\time-1)} }\Pr\left(e \notin G_{\dvec}\mid 
					F\subseteq G_{\dvec}\subseteq H^{(\time-1)}  \right)} & \text{otherwise.}
			\end{cases}
		\]
		Note that if $H^{(\time-1)}$ has  more than $M$ edges
		then 
		$\Pr\left(F \subseteq G_{\dvec} \subseteq  H^{(\time-1)} \right)>0$
		and there is an edge $e \in H^{(\time-1)} $
		such that $\Pr\left(e \notin G_{\dvec} \mid 
		F\subseteq G_{\dvec}\subseteq  H^{(\time-1)}\right)>0$.
		Therefore $H^{(\time)}$ eventually becomes a graph with degree
		sequence $\dvec$, which we denote by $H^{(\infty)}$. 

   Similarly to Claim~\ref{claim:crucial}, we describe the distribution of  $H^{(\tau)}$ at the steps where the number of edges is decreased.  
    For   $i\leq N-M$, we consider  model $\G^+_F (n,\dvec, i)$,   in which  we take a graph 
	$G_{\dvec} \sim \G(n,\dvec)$ 
    conditioned to the event that $F\subseteq G_{\dvec}$ and increase its edge set 
    by    $N-M-i$ edges chosen uniformly at random from the complement of $G_{\dvec}$. Equivalently, a graph from $\G^+_F (n,\dvec, i)$ is the compliment of a random set of $i$ edges sampled as follows:  take a uniform random element  $G_{\dvec,F}$ of the set of graphs with degree sequence $\dvec$ containing $F$ and then pick $i$ edges uniformly from  $K_n \setminus G_{\dvec,F}$.
   % Let $\dvec^{\,c}$ denote the complementary degree sequence to $\dvec$, that is, all components of $\dvec+ \dvec^{\,c}$ equal $n-1$.
   \begin{claim}\label{claim:crucial2}
    If $ i \in [N-M]$ and  $\tau_i$ is the first step where 
    $K_n \setminus H^{(\tau)}$ has $i$ edges  then $H^{(\tau_i)} $  is distributed according 
    $\G_F^+(n,\dvec,i)$.
   \end{claim}
   We prove Claim \ref{claim:crucial2} together with Claim \ref{claim:crucial} in Section \ref{S:crucial}.  
		 For now, we proceed with the proof of Theorem \ref{T:top} assuming that Claim \ref{claim:crucial2} is true.

First, we establish some useful bounds about the quantities $ \hat{\xi}$ and $\zeta$ defined by   
\begin{equation*}%\label{def:hatxi}
			 \hat{\xi} := (\xi \log n)^{1/3} \dfrac{d}{n} \qquad \text{and} \qquad \zeta:= (\xi/\log n)^{1/2},
		\end{equation*}
which we will repeatedly use in the proof.
By the assumptions, we derive that
  		\begin{equation}\label{ineq:hatxi}
			%\zeta^{1/2} \hat{\xi} 
 \log^3 n \ll \dfrac{\xi d  }{\zeta} \ll \hat{\xi}n \ll  d. 
		\end{equation}
The first bound in \eqref{ineq:hatxi} is equivalent to 
$\xi d^2 \gg \log^5 n$ which is true by the assumptions 
$d\gg \xi^{-3} \log n$ and $\xi = o(1/\log n)$.
The  other two bounds are even more straightforward.     Observe  also that $\zeta \ll 1/\log n$ and 
\begin{equation}\label{eq:zeta-log-xi}
   \zeta \log(1/\hat{\xi}) = \zeta \left( \log \left(\dfrac{n }{d\log^{1/3} n} \right)+ \dfrac13 \log(1/\xi) \right) \leq \zeta   \left( \log n+ \log n \right) \ll 1. 
\end{equation}
Now we are ready to proceed. Let   $(G_*, F_*,G_0)$ be the random triple  given by Lemma \ref{L:top-partial}.
Let
		\[ % begin{equation} \label{def:mu}
			\Time \sim \po(\mu), \qquad  \mu:= - N \log {\hat{\xi}}.
		\]  
		We will show that the coupling  $(G, G^*):= (H^{(\infty)}, H_{\zeta}^{(\Time)}  \cup G_0)$  satisfies the theorem requirements, where  
		$H^{(\infty)}$ and $H_{\zeta}^{(\Time)}$ are defined above with $F = F_*$. 
         
        Recall from Lemma \ref{L:top-partial} that $F_*\sim \G(n,\dvec,m_*)$ for some  $m_* = (1-O(\xi))M$.
		Using Claim~\ref{claim:crucial2}, for any    $G'$ with degree sequence $\dvec$,  we get that
		\begin{align*}
			\Pr(H^{(\infty)} = G') &= \sum_{F \subseteq G'} \Pr(G_{\dvec} = G' \mid F \subseteq G_{\dvec}) \cdot \Pr(F_* = F)
			\\
			&= \sum_{F \subseteq G'} \Pr(G_{\dvec} = G' \mid F \subseteq G_{\dvec}) \cdot 
			\frac{\Pr(F \subseteq  G_{\dvec})}{\binom{M}{m_*^{}}}  
       \\
       &= \frac{1}{\binom{M}{m_*^{}}}\sum_{F \subseteq G'} \Pr(G_{\dvec} = G') = \Pr(G_{\dvec} = G').
			% \\&=  \Pr(G_{\dvec} = G') \cdot \sum_{F' \subseteq G'}   \binom{M}{|F'|}^{-1}  \Pr(|F|=|F'|)
			% \\
			% &=\Pr(G_{\dvec} = G') \cdot \sum_{m'=0}^M \sum_{\substack{F'\subseteq G'\\ |F'|=m'}}  \binom{M}{|F'|}^{-1} \Pr(|F|=m') \\
			% &=  \Pr(G_{\dvec} = G') \cdot  \sum_{m'=0}^M   \Pr(|F|=m') = \Pr(G_{\dvec} = G').
		\end{align*}
		Therefore $H^{(\infty)} \sim \G(n,\dvec)$ as required.
		
		Next, note that the construction of the sequences  $K_n \setminus H_0^{(\time)}$ and $K_n \setminus H_\zeta^{(\time)}$
		is identical to that for the graphs $G_0^{(\time)}$ and $G_\zeta^{(\time)}$ 
		considered in the proof of  Theorem \ref{thm:bottom} (but in this proof we use a different value of $\zeta$).  Let
		\[ \hat{m}^{(\time)} = |H^{(\time)}|, \qquad 
		\hat{m}_0^{(\time)} = |H_0^{(\time)}|,
		\qquad 	\hat{m}_\zeta^{(\time)} = |H_\zeta^{(\time)}|. 
		\]
		For convenience, we rewrite the bounds  from \cite{gao2020Sandwiching}
		in these notations.
		By  \cite[Lemma~7]{gao2020Sandwiching},  
		\begin{equation}\label{dist:H}
			H_0^{(\time)}\sim \G(n,\hat{m}_0^{(\time)}), \qquad  H_\zeta^{(\time)}\sim \G(n,\hat{m}_\zeta^{(\time)})
		\end{equation}
		and,  provided  $\time \zeta = o(N)$, we have that
		\begin{equation}\label{eq:m-mH}
			\hat{m}_0^{(\time)}   = (1+O(\zeta + \time \zeta/N)) \hat{m}_\zeta^{(\time)} 
		\end{equation}
		with probability at least 
		$1 - e^{-\Omega \left( \zeta^2 (N-M)^2/\time\right)}$. 
        From \eqref{ineq:hatxi} and the assumption that $d\gg \xi^{-3}\log n$, we get  $\hat{\xi} = o(1)$ and $\hat \xi n  \geq \xi d \gg 1$. Thus, we can estimate
		\begin{equation}\label{range-mu}
			\mu =  N \log (1/\hat \xi)  \in [N, N \log n].
		\end{equation}
		In particular, since $\zeta \ll 1/\log n$,  we get $e^{-\Omega \left( \zeta^2 (N-M)^2/\time\right)} = e^{-\omega (\log n)}$ for $\time = O(\mu)$.
        Using Lemma \ref{l:technical}  and \eqref{range-mu}, we estimate that 
		\begin{equation}\label{Tau-estimate}
		\Pr(\Time > 2\mu ) = e^{-\Omega(\mu)}  =e^{-\omega(\log n )}.
		\end{equation}
        Therefore, we get that \eqref{eq:m-mH} holds for all steps $\tau \leq \Time$ with probability at least $1- e^{-\omega(\log n)}$

	Next, by  \cite[Lemma~2]{gao2020Sandwiching},  we have
		\[
		H_0^{(\Time)}  \sim \G(n,\hat{\xi}), \qquad H_\zeta^{(\Time)}  \sim 
		\G(n,\hat{p}_{\zeta}).
		\] 
		where $\hat{p}_{\zeta}:= \hat{\xi}^{1-\zeta}.$  From  \eqref{eq:zeta-log-xi}, we get 
		$\mu  \zeta = N \zeta \,	{\log (1/\hat{\xi})} = o(N)$. Using also \eqref{ineq:hatxi}, we find that
		\[
		\hat{p}_{\zeta} =\hat\xi \exp\left(\zeta\log (1/\hat\xi)\right)  = (1+o(1)) \hat{\xi}  = o\left(\dfrac{d}{n}\right).
		\] 
		Note that $H_{\zeta}^{(\Time)}$
		is constructed independently of $G_0 \sim \calG(n,p_0)$. Therefore,  we get that 
		$H_{\zeta}^{(\Time)} \cup G_0 \sim \calG(n,p^*)$, where
		\[
		p^* = 1 - (1-p_0)(1-\hat{p}_{\zeta})  = \left(1+ O(\xi)\right) \dfrac{d}{n}+  o\left(\dfrac{d}{n}\right)= (1+o(1)) \dfrac{d}{n}.
		\]
		as required.

		It remains to prove  the claimed probability bound for the event 
		$H^{(\infty)}\not\subseteq H_{\zeta}^{(\Time)} \cup G_0$. 
		Applying the Chernoff bound for $|\G(n,\hat{\xi})|\sim \Bin(N,\hat{\xi})$ (see Lemma \ref{l:technical}), we get that  
		\begin{equation}\label{eq:lowM'}
			\hat m^{(\Time)}, \hat m_{\zeta}^{(\Time)} \geq \hat m_0^{(\Time)} \sim |\G(n,\hat{\xi})| \geq \hat{\xi} N/2
		\end{equation}
		with probability at least $1- e^{-\Omega(\hat{\xi} N)}\geq 1-e^{-\omega(\log n)}$.
		% Using Lemma \ref{l:technical} again and \eqref{range-mu}, we estimate that 
		% \[
		% \Pr(\Time > 2\mu ) = e^{-\Omega(\mu)}  =e^{-\omega(\log n )}.
		% \]
 Note that \eqref{eq:lowM'} implies 
 \begin{equation}\label{tau_i-Tau}
    \tau_i \leq \Time \text{ for all $i\leq (1-\hat{\xi}/2)N -M$.}
        \end{equation}
		The required probability is estimated  similarly to \eqref{eq:union:sparse}:
		\begin{align*} 
				\Pr(H^{(\infty)} \not\subseteq{}& H_{\zeta}^{(\Time)} \cup G_0) \leq \Pr(H^{(\Time)} \not\subseteq H_{\zeta}^{(\Time)} \cup F_*) \\
				&\leq   e^{-\omega(\log n)} +   \sum_{ i\leq (1-\hat{\xi}/2)N -M} \, \Pr\left(\max_{e\in H^{(\tau_i)}\setminus F_* }\left\{\eta^{}_{F_*}(H^{(\tau_i)},e)\right\}>\zeta,\ H^{(\time_i)} \subseteq H_{\zeta}^{(\time_i)} \cup F_*  \right).
		\end{align*}
		It is sufficient to show that  for all $i\leq (1-\hat{\xi}/2)N -M$, we have  
	\begin{equation}\label{top-sufficient}\Pr\left(\max_{e\in H^{(\tau_i)}\setminus F_* }\left\{\eta^{}_{F_*}(H^{(\time_i)},e)\right\}>\zeta,   \
		H^{(\time_i)} \subseteq H_{\zeta}^{(\time_i)} \cup F_* \right) \leq e^{-\omega(\log n)}.
        \end{equation}
		We postpone the proof of the following claim to Section~\ref{sec:proof_of_co-sparse} as it is only about applying Theorem \ref{l:co-sparse:simple} and verifying its assumptions.

		\begin{claim}
  \label{Claim:co-sparse}

%\jc{Is the following easier:}
Let $i\leq (1-\hat{\xi}/2)N -M$. With probability at least $1-e^{-\omega(\log n)}$, the following implication holds:
if $H^{(\time_i)} \subseteq H_{\zeta}^{(\time_i)} \cup F_*$ then
			\begin{equation*}
			\Pr\left(e \in G_{\dvec} \mid 
			F_*\subseteq G_{\dvec} \subseteq H^{(\time_i)} \right) = 
			(1+o(1)) \dfrac{M-m_*}{N-i - m_*}
			\end{equation*}
uniformly for all   $e \in  H^{(\time_i)} \setminus F_*$.

        \end{claim}

We proceed to the final part of the proof of Theorem \ref{T:top}.
Let $\mathcal{E}^{(\time_i)}$ denote the event 
appearing in the implication in Claim \ref{Claim:co-sparse} that 
			\begin{equation*}
			\Pr\left(e \in G_{\dvec} \mid 
			F_*\subseteq G_{\dvec} \subseteq H^{(\time_i)} \right) = 
			(1+o(1)) \dfrac{M-m_*}{ N- i-m_*},\quad \text{for all  $e \in  H^{(\time_i)} \setminus F_*$.}  \end{equation*}
Notice that $\mathcal{E}^{(\time_i)}$ implies that $\max_{e\in \in H^{(\time_i)}\setminus F_*} \eta^{}_{F_*}(H^{(\time_i)},e)\leq \zeta$. Indeed, 
using   the middle inequality of \eqref{ineq:hatxi} and inequality 
$N-i - m_*\geq \hat{\xi}N/2 +M-m_* \geq \hat{\xi}N/2$ and
recalling from   Lemma \ref{L:top-partial} that
$M-m_* = O(\xi M)$, we get
\[
\eta_{F_*}^{}(H^{(\time_i)},e)=1-\frac{1-(1+o(1)) \dfrac{M-m_*}{N-i-m_*}}{1-(1+o(1)) \dfrac{M-m_*}{N-i-m_*}} = o\left(\dfrac{M-m_*}{N-i-m_*}\right)
  = o\left(\dfrac{\xi M}{\hat \xi N}\right)  
  =o\left(\dfrac{\xi d}{\hat \xi n}\right) 
  \leq \zeta 
\]
for all  $e\in  H^{(\time_i)}\setminus F_*$.
Therefore, with $\overline{\mathcal{E}^{(\time_i)}}$ denoting the complement of $\mathcal{E}^{(\time_i)}$,
\begin{align*}
&\Pr\left(\max_{e\in H^{(\tau_i)}\setminus F_* }\left\{\eta^{}_{F_*}(H^{(\time_i)},e)\right\}>\zeta,   \
		H^{(\time_i)} \subseteq H_{\zeta}^{(\time_i)} \cup F_* \right)
        \\
        &\hspace{2cm}\leq \Pr\Big( \overline{ \mathcal{E}^{(\time_i)}}, H^{(\time_i)} \subseteq H_{\zeta}^{(\time_i)} \cup F_* 
        \Big) \\
        &\hspace{2cm}= \Pr\Big(
		\text{the implication } (H^{(\time_i)} \subseteq H_{\zeta}^{(\time_i)} \cup F_*) 
        \longrightarrow \mathcal{E}^{(\time_i)} \text{ fails}\Big).
\end{align*} 
Using  Claim \ref{Claim:co-sparse}, we get \eqref{top-sufficient}, completing the proof.
\end{proof}

\subsection{Proof of Claim~\ref{claim:crucial} and Claim~\ref{claim:crucial2}}\label{S:crucial}
First, we observe that a special case of Claim~\ref{claim:crucial2} with  $F = \emptyset$ reduces to Claim~\ref{claim:crucial} for a complementary degree sequence.  
Recall that the complementary degree sequence  $\dvec^{\,c}$  is such that all components of  $\dvec+\dvec^{\,c}$ are equal to $n-1$.
Comparing the two constructions, one can see that if  $(H^{(\time_i)})_{i=0,\ldots,N-M}$ is a sequence defined as 
in Section \ref{s:top1}  for $F=\emptyset$  
then the sequence of its complements  $(K_n\setminus H^{(\time_i)})_{i = 0\ldots N-M}$ is equivalent to the sequence $(G^{(\time_i)})_{i = 0\ldots N-M}$   defined as in Section \ref{s:bottom} for $\dvec^{\,c}$   instead of $\dvec$. Furthermore, Claim~\ref{claim:crucial2} follows immediately from Lemma~\ref{lem:distribution-H} and Lemma ~\ref{lem:distribution-H2} given in this section.

We consider a relatively simpler process that generates the sequence of graphs in the statement of the lemma. Recall that $M$ is the number of edges in a graph with degree sequence $\dvec$. Let $(H_i)_{0\leq i\leq N-M}$ be defined as follows: $H_0$ is the complete graph on $[n]$, and for every $1\leq i\leq N-M$, for every $e\in H_{i-1}\setminus F_0$, let $H_i=H_{i-1}\setminus e$ with probability 
\begin{equation}
\frac{\Pr(e\notin G_{\dvec} \mid  F\subseteq G_{\dvec} \subseteq H_{i-1} )} {\sum_{e\in H_{i-1}\setminus F}\Pr(e\notin G_{\dvec} \mid  F\subseteq G_{\dvec} \subseteq H_{i-1})}.\label{eq:edgeProb}
\end{equation}
The next lemma established the relation of 
$(H_i)_{0\leq i\leq N-M}$ with the sequence of random graphs considered in Section \ref{s:top1}.
 
\begin{lemma}\label{lem:distribution-H2}
    The sequence $(H^{(\tau_i)})_{0\leq i\leq N-M}$ has the same distribution as $(H_i)_{0\leq i\leq N-M}$.
\end{lemma}
 \begin{proof}
 Recall that for every $0\leq i\leq N-M$,  $\tau_i$ is the minimum $\time$ such that $|K_n\setminus H^{(\time)}|=i$. 
 By our construction,   for any integer $j\geq 1$,  any graph $H$ containing $F$ 
 with $N-i+1$ edges, and any $e \in H \setminus F$, 
\begin{align*}
    \Pr(H^{(\tau+j)}&=H\setminus e \mid H^{(\tau)}=H,\  \tau_{i-1}=\tau, \ \tau_{i}=\tau+j) 
    \\ &= 
     \Pr\left(e^{(\tau+j)} = e \mid a^{(\tau+j)} < 1-\eta_F(H,e^{(\tau+j)}),\, H^{(\tau)}=H,\  \tau_{i-1}=\tau, \ \tau_{i}>\tau+j-1\right)
     \\ &=\Pr\left(e^{(\tau+j)} = e \mid a^{(\tau +j)} < 1- \eta_F(H,e^{(\tau+j)})\right)\\
     & = \frac{\Pr\left(a^{(\tau+j)} < 1- \eta_F(H,e)\mid e^{(\tau+j)} = e \right) \cdot \Pr\left(e^{(\tau+j)} = e \right)}{
     \Pr\left( a^{(\tau+j)} < 1- \eta_F(H,e^{(\tau+j)})\right)}
     \\ & =  \frac{(1- \eta_F(H,e)) \cdot \frac{1}{N}}{\Pr\left( a^{(\tau+j)} < 1- \eta_F(H,e^{(\tau+j)})\right)} \\&\propto
     \Pr(e\notin G_{\dvec} \mid  F\subseteq G_{\dvec} \subseteq H ),\quad \text{(by definition of $\eta_F(H,e)$).}
\end{align*}
In particular, we used that 
   $e^{(\tau + j)}$ is independent of $a^{(\tau+j)}$ and all $e^{(\tau')}$, $a^{(\tau')}$ with $\tau
'< \tau + j$.
Since these conditional probabilities sum up to $1$ over all 
 $e \in H \setminus F$, we get that 
\begin{equation*}
\Pr(H^{(\tau+j)}=H\setminus e \mid H^{(\tau)}=H,\  \tau_{i-1}=\tau, \ \tau_{i}=\tau+j)=\frac{\Pr(e\notin G_{\dvec} \mid  F\subseteq G_{\dvec} \subseteq H )} {\sum_{e\in H\setminus F}\Pr(e\notin G_{\dvec} \mid  F\subseteq G_{\dvec} \subseteq H)}. 
\end{equation*}
Now, it is easy to find the conditional distribution for $H^{(\tau_i)}$ given $H^{(\tau_{i-1})}$:
\begin{align*}
\Pr&(H^{(\tau_i)}=H\setminus e \mid H^{(\tau_{i-1})}=H)  
\\ &= 
\sum_{\tau, j} \Pr(H^{(\tau_i)}=H\setminus e,\, \  \tau_{i-1}=\tau, \ \tau_{i}=\tau+j \mid H^{(\tau_{i-1})}=H)
\\
&=\sum_{\tau,j} \Pr(H^{(\tau+j)}=H\setminus e \mid H^{(\tau)}=H,\   \tau_{i-1}=\tau, \ \tau_{i}=\tau+j)
\\ &\hspace{4cm} \times
\Pr( \tau_{i-1}=\tau, \ \tau_{i}=\tau+j \mid H^{(\tau_{i-1})}=H)
\\
&= \frac{\Pr(e\notin G_{\dvec} \mid  F\subseteq G_{\dvec} \subseteq H )} {\sum_{e\in H\setminus F}\Pr(e\notin G_{\dvec} \mid  F\subseteq G_{\dvec} \subseteq H)}  \sum_{\tau, j}\Pr( \tau_{i-1}=\tau, \ \tau_{i}=\tau+j \mid H^{(\tau_{i-1})}=H)
\\&=\frac{\Pr(e\notin G_{\dvec} \mid  F\subseteq G_{\dvec} \subseteq H )} {\sum_{e\in H\setminus F}\Pr(e\notin G_{\dvec} \mid  F\subseteq G_{\dvec} \subseteq H)}.
\end{align*}
Since this recurrence relation coincides with \eqref{eq:edgeProb}, we conclude that the two sequences have the same distribution.
 \end{proof}

The following lemma basically repeats the counting arguments of \cite[Lemma 3]{gao2020Sandwiching}.
\begin{lemma}\label{lem:distribution-H}
    For every $0\leq i\leq N-M$, we have $H_i \sim \G^+_F(n,\dvec,i)$.
\end{lemma}

\begin{proof}
For each $0\leq i<N-M$,
let $e_{i+1}=H_{i}\setminus H_{i+1}$.
We prove by induction on $i$ that $H_i\sim \G^+_F(n,\dvec,i)$.
The base case $i=0$ is trivially true. Assume the assertion holds for some $i\geq 0$. Given  a graph $H$ whose degree sequence $\degvec(H)\geq \dvec$ and $|K_n\setminus H|=i+1$, let $\S_F(H)$ be the set of graphs $S$ with degree sequence $\dvec$ such that $F\subseteq S\subseteq H$. By definition, to verify that $H_{i+1}\sim \G_F^+(n,\dvec,i+1)$, one needs to show that for every $G$ with $N-i-1$ edges containing $F$, we have 
\[
\Pr(H_{i+1}=G) = \sum_{S\in\S_F(G)} \frac{\Pr(G_{\dvec}=S)}{\binom{N-M}{i+1}}.
\]
 Since $\Pr(G_{\dvec}=S)$ and $\binom{N-M
}{i+1}$ are both independent of $G$, verifying the above equation is equivalent to  
\begin{equation}
\Pr(H_{i+1}=G) \propto |\S_F(G)|.  \label{eq:verify}
\end{equation}

To verify~\eqref{eq:verify}, let ${\mathcal H}(G)$ be the set of supergraphs $H \supseteq G$  with $|H|=|G|+1 = N-i$.  
By the definition of $(H_i)_{0\leq i\leq N-M}$ 
  and the induction hypothesis with (see \eqref{eq:edgeProb} and \eqref{eq:verify}), we find that 
\begin{align*}
\Pr(H_{i+1}=G)&=\sum_{H\in {\mathcal H}(G)} \Pr(H_i=H)\Pr(e_{i+1}=H\setminus G)\\
&\propto\sum_{H\in {\mathcal H}(G)} |\S_F(H)|\; \Pr(e_{i+1}\notin G_{\dvec} \mid F\subseteq G_{\dvec} \subseteq H).
\end{align*}
Next, by definition of $\S_F$,  we get that 
\[
\Pr(e_{i+1}\notin G_{\dvec} \mid F\subseteq G_{\dvec} \subseteq H) = \frac{|\S_F(G)|}{|\S_F(H)|}.
\]
Therefore, 
\[
    \Pr(H_{i+1}=G) \propto 
    \sum_{H\in {\mathcal H}(G)} |\S_F(H)| \cdot \frac{|\S_F(G)|}{|\S_F(H)|} = (N-i-1) |\S_F(G)|\propto |\S_F(G)|,
\]
so we get (\ref{eq:verify}). This verifies the induction step and completes the proof.
\end{proof}

%%%%%%%
%%%%%
%%%%%%%%%%

%%%%%%%%%%%
	%%%%%%%%%%%%%
%%%%%%%%%%%
	%%%%%%%%%%%%%

	%%%%%%%%%%%
	%%%%%%%%%%%%%
	\subsection {Proof of Claim \ref{Claim:co-sparse}}\label{s:top2}
	\label{sec:proof_of_co-sparse}
	Recall that   $(G_*, F_*,G_0)$ is the random triple  given by Lemma \ref{L:top-partial}.
	Define
	\[ % begin{equation}\label{eq:SSt}
		S:=H^{(\time_i)}\setminus F_*,\qquad \tvec := \dvec - \degvec(F_*).
	\]
	Since every subgraph of $H^{(\time_i)}$ with degree sequence $\dvec$ containing $F_*$ uniquely corresponds to a $\tvec$-factor of $S$,  the required probability can be rewritten as follows:
	\[ % begin{equation}\label{eq:Pr_rew}
		\Pr\left(e \in G_{\dvec}\mid 
		F_*\subseteq G_{\dvec} \subseteq H^{(\time_i)}\right) =  \Pr(e \in S_{\tvec} ).  
	\]
	It remains to show that all edge probabilities of random $\tvec$-factor of $S$ are approximately the same  if $H^{(\time_i)} \subseteq H_{\zeta}^{(\time_i)} \cup F_*$ up to a rare event of probability $e^{-\omega(\log n)}$.
	We will do it by applying Theorem \ref{l:co-sparse:simple}.
	The following three lemmas will be useful for verifying its assumptions.

	\begin{lemma}\label{J:lemma}
		Suppose there exist a supergraph $S'\supseteq S$ and  $c\in(0,1)$ such that 
		\[
		\|A(S') - \tfrac{|S'|}{N}J \|_2 + 2\varDelta(S')-2\varDelta(S)+ 2\range{S} \leq c \varDelta(S). 
		\]
		Then $S$ satisfies assumption (ii) of Theorem~\ref{l:co-sparse:simple} with the same $c$. 
	\end{lemma}	

	\begin{proof}
		By the triangle inequality, we get that 
		\[
		\|A(S) - \tfrac{|S|}{N}J \|_2 
		\leq \|A(S') - \tfrac{|S'|}{N}J \|_2 + 
		\|A(S'\setminus S) - \tfrac{|S'\setminus S|}{N}J \|_2
		\]
		Observe that, for any graph $G$, 
		\begin{equation}\label{eq:Gspec}
			\begin{aligned}
				\|A(G) - \tfrac{|G|}{N}J\|_2 
				&\leq \|A(G) - \tfrac{|G|}{N}J\|_1\\
				&= \|(1-\tfrac{|G|}{N})A(G)  -\tfrac{|G|}{N}(I+A(K_n\setminus G))   \|_1
				\\
				&\leq \max\left\{ (1-\tfrac{|G|}{N}) \|A(G)\|_1+ \tfrac{|G|}{N}  (n-1), \tfrac{|G|}{N}  n\right\}
				\\&\leq 2 \varDelta(G). 
			\end{aligned}
		\end{equation}
		In the  penultimate inequality of \eqref{eq:Gspec}, we observed that any non-zero row of $A(G)$ corresponds to a row of $I+A(K_n\setminus G))$ with at most $n-1$ non-zero entries.
		Applying \eqref{eq:Gspec} with $G= S'\setminus S$ we get that 
		\[
		\|A(S'\setminus S) - \tfrac{|S'\setminus S|}{N}J \|_2
		\leq 2 \varDelta(S'\setminus S) \leq 2(\varDelta(S') - \varDelta(S) + \range{S}),
		\]
        which gives  assumption (ii) of Theorem~\ref{l:co-sparse:simple}.
	\end{proof}
	
	\begin{lemma}\label{l:Z_xy}
		Let $H \sim \G(n,m)$ where  $m \gg n\log n$.
		Let $h_1,\ldots,h_n$  be the degrees of the random graph $H$ and 	  let $A(H)$ denote  its  adjacency matrix.    
		Assume 
		$\eps = \eps(n)>0$ is such that
		$ \sqrt{\dfrac{n\log n}{m}}\ll \eps \ll 1$.
		Then, with probability $1 -  e^{-\Omega( \eps^2 m/n)}$,  
		\[
		\max_j |h_j  - 2m/n| \leq \eps m/n, \qquad 
		\|A(H) - (m/N) J\|_2 \leq \eps m/n,
		\]
		and, 
		uniformly	for all $\xvec,\yvec \in [0,1]^n$  with $\|\xvec\|_1, \|\yvec\|_1 = \Omega(n)$,
		\[
		\biggl | \sum_{(jk)\st jk \in H} x_j y_k- (m/N)\|\xvec\|_1\|\yvec\|_1 \biggr | \leq \eps m.
		\]
	\end{lemma}
	\begin{proof}
		It is sufficient to prove the bounds above for  the random graph $\tilde{H}\sim \G(n,p)$  where $p=m/N$, 
		since the probability of the  event  $|E(\tilde{H})|=m$ is 
		substantially larger than  $e^{-\Omega(\eps^2 m/n)}$. 
		Observe that the degrees of  $\tilde{H}$ are distributed according to $\Bin(n-1,p)$.
		Applying Lemma~\ref{l:technical} and using the union bound, we
		obtain the uniform concentration bound for $h_j-2m/n$. 
		The two other bounds  for $\tilde{H}$ hold with even better probability estimates 
		and are given in  Lemma \ref{A:spec} and Lemma \ref{A:z_con} in the Appendix.
	\end{proof}

	\begin{lemma}\label{lem:union}
		Let $G_1\sim\G(n,m)$ and $G_2\sim\G(n,p)$ where $G_1$ and $G_2$ are independent. Let $G_{\ominus}=G_1\setminus G_{2}$  then $G_{\ominus} \sim \G(n,|G_{\ominus}|)$ and 
		$|G_{\ominus}| \sim \Bin(m,1-p)$.
	\end{lemma}
	
	\begin{proof}
		For any graph $G'$ with $m' \leq m$ edges we have
		\begin{align*}
			\Pr(G_{\ominus} = G')&= \sum_{\substack {G_1' \supseteq G'\\ |G_1'|=m}}  \Pr(G_1  \setminus G_1' = G' ) \cdot \Pr(G_1 = G_1')
			\\
			&= \binom{N-m'}{m-m'} \cdot p^{m'}(1-p)^{m-m'} \cdot \binom{N}{m}^{\!-1}.
		\end{align*}
		Clearly, this probability does not depend on the structure of $G'$ but only depends on the number of its edges. Multiplying by $\binom{m}{m'}$, which is the number of graphs $G'$ with $m'$ edges, we get the desired distribution of $|G_{\ominus}| \sim \Bin(m,1-p)$.
	\end{proof}

	We proceed to the proof of  Claim~\ref{Claim:co-sparse}. 
	Let 
	\[
	S':=  H_{\zeta}^{(\time_i)}, \qquad  S_0:= H_0^{(\time_i)} \setminus G_0.
	\]
	Note that if 
	$H^{(\time_i)} \subseteq H_{\zeta}^{(\time_i)} \cup F_*$ then
	\begin{equation}\label{recall:contain}
		S_0  \subseteq  S \subseteq S'. 
	\end{equation}
	Indeed, the upper containment is immediate and the lower containment always holds by our construction since $H_0^{(\time_i)} \subseteq S\setminus F_*$ and  $F_* \subseteq G_0$ by Lemma \ref{L:top-partial}(a). Recall  that  $G_0 \sim \G(n,p_0)$ with $p_0 = (1+o(1))d/n$.
	Using \eqref{dist:H} and  Lemma~\ref{lem:union}, we get that  
	\begin{equation}
	S'\sim \G(n,\hat{m}_{\zeta}^{(\time_i)}), \qquad S_0\sim \G(n, |S_0|), \qquad   |S_0| \sim \Bin(\hat{m}_0^{(\time_i)}, 1-p_0),\label{eq:distr-S-S0}
	\end{equation}
 where $\Bin(X,q)$ for a positive integer-valued random variable $X$ denotes the compound  binomial distribution of $\Bin(x,q)$ with parameter $x$ distributed as $X$.
	From \eqref{eq:lowM'}, \eqref{tau_i-Tau}, and the fact that $\hat{m}_{\zeta}^{(\time)}$ and $\hat{m}_{0}^{(\time)}$ are decreasing functions of $\time$, we have that 
	\begin{equation}
	\hat{m}_{\zeta}^{(\time_i)}\geq \hat{m}_{0}^{(\time_i)} \geq \hat{\xi}N/2 \label{eq:hat-m-bound}
	\end{equation}
	with probability at least $1-e^{-\omega(\log n)}$. 
    Using \eqref{ineq:hatxi}, we easily get
\[
    \hat{m}_{0}^{(\time_i)} p_0  = \Omega( \hat{\xi} N  d/n)  \gg \log n.
\]
 From  \eqref{eq:zeta-log-xi}, \eqref{Tau-estimate}, \eqref{tau_i-Tau},  we get  
		\[
        \time_i \zeta/N \leq 2\mu \zeta /N = o(1).
        \]
        with the same probability bound.
     Using \eqref{eq:m-mH} and the Chernoff bound  (see  Lemma~\ref{l:technical}), we obtain that 
	\begin{equation*}%\label{diff_SS}
		|S_0| =(1+O(\zeta + \time_i \zeta/N)) |S'| =  (1+o(1))|S'|
	\end{equation*}
	with probability at least $1- e^{-\omega(\log n)}$.  Due to  \eqref{recall:contain},
 with probability at least $1- e^{-\omega(\log n)}$,
 	\begin{equation}\label{diff_SS}
		|S|,|S_0| 
        %=(1+O(\zeta + \time \zeta/N)) |S'| 
        =  
        (1+o(1))|S'|.
	\end{equation}

	Let $\svec = (s_1, \ldots, s_n)$ be the degree sequence of $S$.  It follows from Claim \ref{claim:crucial2} and the definition of $\calG_F(n,\dvec,i)$ that 
	\begin{equation}\label{dist_s}
		s_j -t_j  \sim  \operatorname{Hypergeometric}(N-M, n-1-d_j, N-M-i).
	\end{equation}
	In particular, we find that
	\[
	\E\left[s_j  \mid F_*  \right]=  t_j+ \dfrac{(n-1-d_j)(N-M-i)}{N-M}.
	\]
	From Lemma \ref{L:top-partial}(b), we have that 
    \begin{equation}\label{t-assymp} t_j  = \Theta(\xi d) \qquad \text{and} \qquad \range{\tvec} = O(\xi^2 d)
    \end{equation}
    with probability at least $1-e^{-\omega(\log n)}$.

    Now we are ready to show that the assumptions of Theorem \ref{l:co-sparse:simple}   hold with probability
    at least $1-e^{-\omega(\log n)}$. By assumptions we have that  
    $N-i-M\geq \hat{\xi}N/2$.
	Using \eqref{dist_s}, \eqref{t-assymp} and applying the Chernoff bound (see Lemma \ref{l:technical}), we find that
	\[  % begin{equation}\label{range-S}
		\range{S} = O(\xi)	\varDelta(S)  \ll \varDelta(S)/\log n
	\]
	with the same probability bound. Thus,  we get condition (i) of Theorem \ref{l:co-sparse:simple}.

	Recall that $\xi d \gg \xi^3 d \gg \log n$.   
    Combining \eqref{ineq:hatxi}, \eqref{eq:hat-m-bound}, and \eqref{diff_SS}, we find that, with probability  at least $1-e^{-\omega(\log n)}$,
	\[
	\varDelta(\tvec) \varDelta(S)\geq  \varDelta(S) \geq \dfrac{2|S|}{n} = \Omega(\hat\xi n)
	\gg \log^2 n
	\]
	 and
	\begin{align*}
		\dfrac{\varDelta(\tvec)}{\varDelta(S)} = 
		O\bigl(\xi d/\hat{\xi}n\bigr)  \ll \zeta \ll \dfrac{1}{\log n},
	\end{align*}
	Using also \eqref{t-assymp}, we verify conditions (iv) of Theorem \ref{l:co-sparse:simple}.

	Next, applying Lemma \ref{l:Z_xy} with $H = S'$
    and $\eps = \dfrac{1}{\log n}$, we find that 
	with a probability at least 
	\[
	1 - e^{-\omega(\log n)}
    - e^{-\Omega(\eps^2  \hat{\xi} n)}
    = 1 - e^{-\omega(\log n)}
	\] 
  the following bounds hold:
	\[
		\varDelta(S') = (1 + O(\eps))\tfrac{2|S'|}{n} , \qquad 
		\|A(S') - \tfrac{|S'|}{N} J\|_2 =O(\eps) \varDelta(S'),
	\]
	and   uniformly	for all $\xvec,\yvec \in [0,1]^n$  with $\|\xvec\|_1, \|\yvec\|_1 = \Omega(n)$,
	\begin{equation}\label{scalar-S'}
		\sum_{(jk)\st jk \in S'} x_j y_k = (1+O(\eps)) 	\tfrac{|S'|}{N} \|\xvec\|_1\|\yvec\|_1.
	\end{equation}
Using Lemma \ref{J:lemma} and  \eqref{diff_SS}, we get the conditions (ii)  of Theorem \ref{l:co-sparse:simple} hold for any constant $c$, for example we can take $c=1/2$.
    
Finally,  by ~(\ref{eq:distr-S-S0}) and Lemma~\ref{l:Z_xy} with $H = S_0$ and $H=S'$ respectively (as done for $H=S'$ in \eqref{scalar-S'}), for all $\xvec,\yvec \in [0,1]^n$  with $\|\xvec\|_1, \|\yvec\|_1 = \Omega(n)$,
\[
	\frac{\sum_{(jk)\st jk \in S_0} x_j y_k}{\|\xvec\|_1 \|\yvec\|_1}= (1+o(1))\dfrac{|S_0|}{N};\quad \text{and}\quad 	\frac{\sum_{(jk)\st jk \in S'} x_j y_k}{\|\xvec\|_1 \|\yvec\|_1}= (1+o(1))\dfrac{|S'|}{N}.
\]
By~(\ref{diff_SS}), it follows that
\[
\frac{\sum_{(jk)\st jk \in S_0} x_j y_k}{\|\xvec\|_1 \|\yvec\|_1}=(1+o(1))\frac{\sum_{(jk)\st jk \in S'} x_j y_k}{\|\xvec\|_1 \|\yvec\|_1}=(1+o(1))\dfrac{|S'|}{N}.
\]
On the other hand, by \eqref{recall:contain}, \[
	\sum_{(jk)\st jk \in S_0} x_j y_k \leq 
	\sum_{(jk)\st jk \in S} x_j y_k \leq \sum_{(jk)\st jk \in S'} x_j y_k.
	\]
Consequently,
	\[
	\frac{\sum_{(jk)\st jk \in S} x_j y_k}{\|\xvec\|_1 \|\yvec\|_1} = (1+o(1)) \dfrac{|S'|}{N}=
	(1+o(1)) \dfrac{|S|}{N} = (1+o(1))\dfrac{\varDelta(S)}{n}.
	\]
	This verifies condition (iii) of  Theorem \ref{l:co-sparse:simple}.

	In the above we checked that all assumptions of Theorem~\ref{l:co-sparse:simple} hold with probability at least $1 - e^{-\omega(\log n)}$.  Then, applying   Theorem~\ref{l:co-sparse:simple} and using~\eqref{diff_SS}, we conclude that 
	\[
	\Pr(e \in S_{\tvec}) =(1+o(1)) \dfrac{\varDelta(\tvec)}{\varDelta(S)} =
       1+o(1)) \dfrac{\|\tvec\|_1}{2|S|}
        =
    (1+o(1))  \dfrac{M-m_*}{N-i-m_*}.
	\]
	This completes the proof of Claim \ref{Claim:co-sparse}.
	
	% % % % % % % % % % 
	% % % % % % % % % %% % % % % % % % % %
	% % % % % % % % % %
	
	\subsection{Proof of Theorem \ref{thm:sandwich2}}\label{s:proof_sandwich}
	The  case $\min\{d,n-d\} = \Theta(n)$ is covered by \cite[Theorem 5]{gao2020Sandwiching}. We can restrict ourselves to the case   $d = o(n)$ since the case $n-d=o(n)$ follows from taking complements.

	Next, we  consider the case when  
	\[ 
	 \log^4 n 
	\ll d \ll \dfrac{n \log \log n}{\log n}.
	\]
	Since $\range{\dvec} \ll \dfrac{d }{\log n}$, there is $\xi = \xi(n)$ such that 
	\[
	\xi \ll \dfrac{ 1}{\log n}, \qquad  \xi n \geq d \gg \xi^{-3} \log n, \qquad \range{\dvec} \leq \xi d.
	\]
	Applying Theorem \ref{T:top} and Theorem \ref{thm:bottom}, we construct two couplings $\pi = (G_*,G)$ and $\pi' = (G,G^*)$ such that $G\sim \G(n,\dvec)$,  $G_* \sim \G(n,p_*)$, $G^*\sim \G(n,p^*)$ with $p_*,p^* = (1+o(1))\dfrac{d}{n}$ and 
	\[
	\Pr(G_* \subseteq G) = 1-o(1), \qquad  \Pr(G \subseteq G^*) = 1-o(1).
	\]
	We can now stitch $\pi$ and $\pi'$
	together to construct the required triple $(G_* , G, G^* )$. First uniformly generate a graph
	$G \sim \G(n,\dvec)$. Then, conditional on $G$, we 
	generate $G_*$ under $\pi$ and generate $G^*$ 
	under $\pi'$.
	
	For the remaining case when
	\begin{equation}\label{case:rem}
		\Omega\left(\dfrac{n \log \log n}{\log n}\right)= d  \ll n,
	\end{equation}
	the bottom side of the sandwich follows from Theorem~\ref{thm:low} or Theorem \ref{thm:bottom}. For the top side,
	we use the embedding from 
	\cite[Theorem 8]{gao2020Sandwiching}, which stated below as Theorem \ref{thm:coupling} for the reader's convenience.
	\begin{theorem}\label{thm:coupling}
		Let $\dvec' = \dvec'(n) \in \Naturals^n$ be a  degree sequence and 
		$\xi = \xi(n)>0$ be such that $\xi(n) = o(1)$.  Denote $\varDelta = \varDelta(\dvec')$. Assume
		$
		n\cdot \range{\dvec'
		} \leq \xi \varDelta(n-\varDelta)$
		and
		$n-\varDelta  \gg \xi \varDelta \gg n / \log n$.
		Then there exist $p = \left(1- O\left(\xi\right)\right)\varDelta/n$, and a coupling  $(G_L,G')$  with
		$G_L \sim \G(n,p)$ and  $G'\sim \G(n,\dvec')$
		such that
		\[
		\Pr(G_L \subseteq G') = 1 - e^{-\Omega \left(\xi^3\varDelta\right)} =1- e^{-\omega(n/\log^3 n) }. 
		\] 
	\end{theorem}
	
	We apply Theorem \ref{thm:coupling} for the complement degree sequence $\dvec'=(n-1)\boldsymbol{1} -\dvec$. Due to \eqref{case:rem} and $\range{\dvec} \ll \dfrac{d }{\log n}$, there is $\xi' = \xi'(n)=o(1)$ such that
	\[
	d \gg \xi' n  \gg \dfrac{n}{\log n}, \qquad \range{\dvec} \leq \xi' d.
	\]
	Observing that $n-\varDelta(\dvec')= (1+o(1))d$ and
	\[
		\varDelta(\dvec') =   n - (1+O(\xi'))d,  
	\]
    we can check that all assumptions of Theorem \ref{thm:coupling}  are satisfied
	with $\xi = \xi'$.
	Let  
	\[
	(G,G^*):= (K_n \setminus G', K_n \setminus G_L),
	\]
	where $(G_L, G')$ is the coupling from Theorem \ref{thm:coupling} where $G_L\sim \G(n,p)$ for some 
	\[
	p=(1-O(\xi'))\dfrac{\varDelta(\dvec')}{n}=(1-O(\xi'))\left(1-(1+O(\xi'))\dfrac{d}{n}\right)=1-(1+O(\xi'))\dfrac{d}{n}+O(\xi').
	\]
	Since $d\gg \xi' n$, it now follows that
	$1-p=(1+o(1)) d/n.
	$    
	Thus, from Theorem \ref{thm:coupling} we get that the random graphs $G,G^*$ have the required marginal distributions and  $G\subseteq G^*$  with probability tending to $1$.
	
	Stitching together the two  couplings from applying Theorem \ref{thm:coupling} as explained in the previous case, we have completed the proof of Theorem \ref{thm:sandwich2}.\qed
	% % % % % % % % % %
	% % % % % % % % % %
	% % % % % % % % % %
	% % % % % % % % % %
	% % % % % % % % % %

%\section{Checking the claims}

	% % % % % % % % % %
	% % % % % % % % % %
	% % % % % % % % % %
	% % % % % % % % % %
	% % % % % % % % % %

\section{Switchings}\label{s:switchings}

In this section we prove the following theorem, which implies Theorem \ref{l:co-sparse:simple}. For $\xvec,\yvec \in \Reals^n$ and a graph $S$ with vertex set $[n]$, denote
	\begin{equation}\label{Z_def}
		\langle  \xvec,\yvec\rangle_S = \sum_{(jk)\st jk \in S} x_j y_k.
	\end{equation}

    	\begin{theorem}\label{l:co-sparse}
		Let $S$ be a graph on $n$ vertices and $\tvec$ be a degree sequence such
		that  the set of $\tvec$-factors of $S$ is not empty and the following assumptions hold for some  $\alpha = \alpha(n) \in (0,1)$ and  some fixed $\beta\geq 1$.
         %$1\leq\beta=\beta(n)=O(1)$
%\jc{To avoid requiring to justify that $\beta\geq 1$. }
   % $\gamma = \gamma(n)=O(1)$, fixed $\beta\geq 1$ and any sufficiently small  $\delta = \delta(\beta)$.
 \begin{itemize}
			\item[\rm (A1)] 
			$\displaystyle \frac{\varDelta(\tvec)}{ \alpha \varDelta(S)} + 
			\frac{1}{ \alpha^2 \varDelta(\tvec) \varDelta(S)} =o(1);
			$
			\item[\rm (A2)]  
			$\displaystyle
			\bigl\|A(S) - \dfrac{|S|}{N}J\bigr \|_2 + 6\range{S}    \leq n^{-\alpha} \varDelta(S);
			$
			\item[\rm (A3)] 
			$	\displaystyle
			\left(\frac{\varDelta(S)}{\varDelta(S)-\range{S}}\right)^{\!1/\alpha} \left(\frac{\varDelta(\tvec)}{\varDelta(\tvec)-\range{\tvec}}\right)^{\!1/\alpha} \leq \beta.
			$
		\end{itemize}
		Then, there is a constant $\delta=\delta(\beta)>0$  such that
        if 
         \begin{equation}\label{def:gamma-S}
    \gamma_S(\delta)= \max_{\substack{\xvec,\yvec \in [0,1]^n\\ \|\xvec\|_1, \|\yvec\|_1\geq \delta n}} \left|\log \frac{\langle\xvec ,\yvec \rangle_S }{\|\xvec\|_1 \|\yvec\|_1 \varDelta(S)/n}\right|<\infty
\end{equation}
        then,
        uniformly over $jk \in S$, 
		\begin{equation}
        \begin{aligned}
		\Pr&(jk \in S_{\tvec}) =\dfrac{\varDelta(\tvec)}{\varDelta(S)} 
        \\ &\times\exp\Bigl(O\Bigl(\gamma_S(\delta)+\dfrac{\varDelta(\tvec)}{\alpha \varDelta(S)}
		+  \dfrac{1}{\alpha^2 \varDelta(\tvec) \varDelta(S)} + \dfrac{\range{S}}{\alpha(\varDelta(S)-\range{S})} +
		\dfrac{\range{\tvec}}{\alpha (\varDelta(\tvec)-\range{\tvec})}\Bigr)\Bigr). \label{eq:edge-prob}
        \end{aligned}
		\end{equation}
	\end{theorem}

	For an edge  $jk \in S$, consider
	the partition of the set of $\tvec$\nobreakdash-factors of $S$ into two disjoint sets $\mathcal{S}(\tvec,jk)$ and  $\mathcal{S}(\tvec, \overline{jk})$, where   elements of  $\mathcal{S}(\tvec,jk)$  contain $jk$ while 
	elements of $\mathcal{S}(\tvec,\overline{jk})$  do not. Since $S_{\tvec}$ is a uniformly random $\tvec$-factor of $S$, we have
	\begin{equation}\label{switch1}
		\Pr (jk \in S_{\tvec}) =  \frac{|\mathcal{S}(\tvec,jk)|}{ |\mathcal{S}(\tvec,jk)| + |\mathcal{S}(\tvec,\overline{jk})|}.
	\end{equation}
	Thus, it is sufficient to estimate the ratio $|\mathcal{S}(\tvec,jk)| / |\mathcal{S}(\tvec,\overline{jk})|$.  We do it  using  the switching method which we briefly describe below.  
	
	We consider  switching operations (corresponding to a certain type of alternating walks) on the set of $\tvec$\nobreakdash-factors of $S$  that replace a $\tvec$-factor containing a given edge $jk$ by another $\tvec$-factor  without   $jk$.
	Given a $\tvec$-factor $T \in \mathcal{S}(\tvec,jk)$, let $\mathcal{F}(T)$ be the set of switching operations that can be applied to $T$ to obtain an element of $\mathcal{S}(\tvec,\overline{jk})$. Potentially there could be multiple switchings  that produce the same $\tvec$-factor of $S$  when applied to $T$.  Similarly, for   $T' \in \mathcal{S}(\tvec,\overline{jk})$ we consider the  set $\mathcal{B}(T')$ of inverse switchings that can be applied to $T'$ to obtain an element of $\mathcal{S}(\tvec,jk)$.  The main idea of  the switching method is  to define the switching operations in such a way that all sets $\mathcal{F}(T)$  
	are of approximately the same size and also all  sets  $\mathcal{B}(T')$ are of approximately the same size.
	Then, using the double counting argument, we can estimate 
	\begin{equation}\label{switch2}
		\frac{\min_{T' \in \mathcal{S}(\tvec,\overline{jk}) } |\mathcal{B}(T')|}{ 
			\max_{T \in \mathcal{S}(\tvec,jk) } |\mathcal{F}(T)|}
		\leq 
		\frac{|\mathcal{S}(\tvec,jk)| }{  |\mathcal{S}(\tvec,\overline{jk})|} 
		\leq \frac{\max_{T' \in \mathcal{S}(\tvec,\overline{jk}) } |\mathcal{B}(T')|}{ 
			\min_{T \in \mathcal{S}(\tvec,jk) } |\mathcal{F}(T)|}.
	\end{equation}

	Next, we formally define our switching operations which are called \textit{$\ell$-switchings}, where $\ell \geq 3$ is an integer.  To perform an $\ell$-switching on a $\tvec$-factor $T\in \mathcal{S}(\tvec,jk)$, choose a sequence of vertices $j=u_1,v_1,u_2,v_2,\ldots,u_{\ell},v_{\ell}=k$ such that
	\begin{itemize} \itemsep=0pt
		\item $u_iv_i$ are edges in $S-T$ for all  $i =1, \ldots, {\ell}$.
		\item $v_iu_{i+1}$ are edges in $T$ for all $i =1, \ldots, {\ell-1}$.
		\item All the above edges are distinct from each other and distinct from $jk$.
		(However, repetition of vertices is permitted.)
	\end{itemize} 
	Then the $\ell$-switching replaces the edges  ${jk}$ and $\{v_iu_{i+1}\}_{i =1, \ldots, {\ell-1}}$  in $T$
	by $\{u_iv_i\}_{i =1, \ldots, {\ell}}$.
	Observe that the resulting graph $T'$ has the same degree sequence as $T$, and  $T' \in \mathcal{S}(\tvec,\overline{jk})$.
	The reverse operation, converting $T'$ to $T$, is called an \textit{inverse $\ell$-switching.}
	See Figure~\ref{f:long} for an illustration.
	\begin{figure}[ht!]
		\centering
		\begin{tikzpicture}[scale=0.6]
			\node (u1) at (0,4) {};
			\node [above left=0 and -0.5 of u1] {$u_1=j$};
			\node (vl)  at (2,4) {};
			\node[above right =0 and -0.5 of vl]{$k=v_{\ell}$};
			\node (v1) [label=left:{$v_1$}]   at (-1.5,2.5) {};
			\node (ul)  [label=right:{$u_{\ell}$}]  at (3.5,2.5) {};
			\node (u2) [label=left:{$u_2$}]  at (-1.5,0.5) {};
			\node (v2) [label=right:{$v_{\ell-1}$}]  at (3.5,0.5) {};
			\node at (1,-1)  {\Huge $\boldsymbol{\ldots}$};
			
			\node at (1,1.6)  {\LARGE $T$};

			\draw [-,thick] (u1) -- (vl);
			\draw [-,dashed](u1) -- (v1);
			\draw[-,thick] (v1) -- (u2);
			\draw[-, thick] (ul)--(v2);
			\draw[-,dashed](ul)--(vl);
			
			\draw [fill] (u1) circle (0.2); \draw [fill] (v1) circle (0.2);  \draw [fill] (u2) circle (0.2);   \draw [fill] (v2) circle (0.2);  \draw [fill] (ul) circle (0.2); 
			\draw [fill] (vl) circle (0.2);
			
			\node   at (5.8,1.5)  {$\Longleftrightarrow$};
			
			\begin{scope}[shift={(9.5,0)}]

				\node (u1) at (0,4) {};
				\node [above left=0 and -0.5 of u1] {$u_1=j$};
				\node (vl)  at (2,4) {};
				\node[above right =0 and -0.5 of vl]{$k=v_{\ell}$};
				\node (v1) [label=left:{$v_1$}]   at (-1.5,2.5) {};
				\node (ul)  [label=right:{$u_{\ell}$}]  at (3.5,2.5) {};
				\node (u2) [label=left:{$u_2$}]  at (-1.5,0.5) {};
				\node (v2) [label=right:{$v_{\ell-1}$}]  at (3.5,0.5) {};
				\node at (1,-1)  {\Huge $\boldsymbol{\ldots}$};
				
				\node at (1,1.6)  {\LARGE $T'$};
				
				\draw[-,dashed] (u1) -- (vl);
				\draw [-,thick] (0,4) -- (-1.5,2.5);
				\draw[-,dashed](v1) -- (u2);
				\draw [-,dashed] (ul)--(v2);
				\draw [-,thick](3.5,2.5)--(2,4);
				
				\draw [fill] (u1) circle (0.2); \draw [fill] (v1) circle (0.2);  \draw [fill] (u2) circle (0.2);   \draw [fill] (v2) circle (0.2);  \draw [fill] (ul) circle (0.2); 
				\draw [fill] (vl) circle (0.2);
				
			\end{scope}
		\end{tikzpicture}
		\caption{$\ell$-switching.}
		\label{f:long}
	\end{figure}

	In the following, $\mathcal{F}_\ell (T)$ is the  set of $\ell$-switchings applicable to 
	a  graph $T\in  \mathcal{S}(\tvec,jk)$. Similarly,  $ \mathcal{B}_\ell(T')$  is the set of inverse $\ell$-switchings applicable to a graph $T'\in \mathcal{S}(\tvec,\overline{jk})$.

	Recall   that $A(G)$ is the adjacency matrix of $G$.	 
	Let $\evec_i$ denote the standard unitary column vector with $1$ in the $i$-th component. 
	For  nonnegative integers $a,b$ define 
	\[
		w_{a,b}(S,T) = \max_{i,P} \|P \evec_i\|_\infty,
	\]
	where the maximum is taken over all $i\in [n]$ and matrices $P$ which are the product of
	$a$ factors $A(S)$ and $b$ factors $A(T)$ (e.g. for $a=1$, $b=2$, the matrix $P$ can be one of 
	$A(S)A(T) A(T)$, $A(T)A(S) A(T)$, $A(T) A(T) A(S)$).   
	Note that the 
	components of $P \evec_i$ 
	correspond to the number of walks that start at $i$ and finish at a given vertex
	which use  $a$ edges from $S$ and $b$ edges from $T$  in a predetermined order (corresponding to $P$). Thus, $w_{a,b}(S,T)$
	is an upper bound  on the number of such walks.

\begin{lemma}\label{l:FB}
	Assume $\ell = 2 h+1$ for $h\geq 2$. 
	Let $A$ be the adjacency matrix of $S$, and define $z=\lfloor(\ell-1)/3\rfloor$. 
	\begin{itemize}
	\item[(a)]     If $T \in \mathcal{S}(\tvec,jk)$,  $B$ is the adjacency matrix of $T$, and    $w_{a,b} = w_{a,b}(S, T)$, then 
		\begin{align*}
			\langle(BA)^{h} \evec_j , (BA)^{h} \evec_k\rangle_S
				&\geq |\mathcal{F}_\ell (T)|   
				\\
				&\geq \langle(BA)^{h} \evec_j , 
				(BA)^{h} \evec_k\rangle_S \\
				&{\quad} - \ell\, w_{\ell-1,\ell} -
				2\ell^2  \varDelta(\tvec)^{\ell-z-2}\varDelta(S)^{\ell-z-1}w_{z,z} .
		\end{align*}
			\item[(b)]   If $T' \in \mathcal{S}(\tvec,\overline{jk})$,  
			$B'$ is the adjacency matrix of $T'$, and    $w_{a,b}' = w_{a,b}(S, T')$, then                   
			\begin{align*}
				\langle(B'A)^{h-1}B' \evec_j , 
				(B'A)^h B'\evec_k\rangle_S
				&\geq 
				|\mathcal{B}_\ell (T')|   
				\\
				&\geq 
				\langle(B'A)^{h-1}B' \evec_j , (B'A)^h B'\evec_k\rangle_S \\
				&{\quad} -	\ell w'_{\ell-2, \ell+1}
				- 2\ell^2  \varDelta(\tvec)^{\ell-z-1}\varDelta(S)^{\ell-z-2}w'_{z,z} .
			\end{align*}
		\end{itemize}
	\end{lemma}
	
\begin{proof}
We start by defining some terminology.
We are considering walks $W =u_1,v_1,\ldots,v_\ell$ in~$S$,
where $u_1=j$ and $v_\ell=k$.
A \textit{labelled vertex} is a vertex in $W$ together with its label.
We refer to labelled vertices by their labels, for example~``$v_2$''.
Since $W$ can repeat vertices, several labelled vertices might be located at
the same vertex of~$S$.
A \textit{subwalk} is a sequence of consecutive labelled vertices of $W$
in either forward or backward order.
If $x$ and $y$ are labelled vertices, then $W[x,y]$ is the subwalk from~$x$
to~$y$.  Also, ``$x=y$'' means that $x$
and $y$ are located at the same vertex of~$S$.

For each $i$ define three subwalks of $2z$ edges:
\[
W_1 = W[u_1,u_{z+1}] ,\quad
W_2(i) = W[u_{i+1},u_{i+z+1}], \quad
W_3 = W[v_{\ell-z},v_\ell].
\]
 
We start with part (a).
Observe that the components of $(BA)^h \evec_j$ 
correspond to counts of walks of length $2h$ which alternate between edges of $S$ and $T$ starting from vertex $j$ and an edge from $S$. 
Clearly, this gives an upper bound for the number of walks that alternate between $S-T$ and~$T$. 
An $\ell$-switching is determined by the subwalks $W[u_1,u_{h+1}]$
and $W[v_\ell,v_{h+1}]$,
plus the edge $u_{h+1}v_{h+1}\in S$.
Summing over all choices of $u_{h+1} v_{h+1} \in S$
and bounding the choices for the two subwalks by the corresponding components
of $(BA)^h \evec_j$
and $(BA)^h \evec_k$,  we have the upper bound for $ |\mathcal{F}_\ell (T)|$.

Some of the walks counted by the upper bound might not be valid $\ell$-switchings.
This could happen under the following conditions:  
\begin{itemize}\itemsep=0pt
   \item[(1)] One of the edges $u_i v_i$ which we choose from $S$ belongs also to  $T$.
   \item [(2)] Condition (1) doesn't hold, but $\{v_i, u_{i+1}\} = \{j,k\}$ for some $i$.
   \item [(3)] Condition (1) doesn't hold, but 
      $\{v_i, u_{i+1}\} = \{v_{i'}, u_{i'+1}\}$ for some $i<i'$ (repeated $T$~edge).       
   \item [(4)] Condition (1) doesn't hold, but 
      $\{u_i, v_i\} = \{u_{i'}, v_{i'}\}$ for some $i<i'$ (repeated $S$~edge).			  
\end{itemize}

Then, by the union bound,
\begin{equation}\label{e:pre-b}
	|\mathcal{F}_\ell (T)|  \geq  \langle(BA)^h \evec_j , 
	(BA)^h \evec_k\rangle_S - N_1 - N_2 -N_3 - N_4,
\end{equation}
where  $N_i$ is the number of invalid choices for $W$
corresponding to condition~($i$).

Consider condition~(1). There are $\ell$ possible choices of $i$ such that 
$\{u_i,v_i\}\in T$, and for each choice $W$ is a walk from~$j$ to~$k$
with $\ell$ edges from $T$ and the remaining $\ell-1$
edges from~$S$, in a particular order.  This implies
\[
    N_1 \leq \ell w_{\ell-1,\ell}.
\]
From now on we will assume that condition~(1) does not hold, so that walks
strictly alternate between $T$ and $S-T$. In particular, walks never
backtrack.

Our arguments for bounding $N_2$, $N_3$ and $N_4$ follow a common pattern.
For each possible edge overlap between $W$ and $jk$ or itself, we define a subwalk $W^*$
of $2z$ edges, and a sequence $\mathcal{Z}$ of subwalks that
(together with the assumed edge overlap)
specify all labelled vertices other than the internal labelled vertices of $W^*$.
By following the subwalks in $\mathcal{Z}$ beginning at $u_1$ or $v_\ell$,
we can choose locations for all of their labelled vertices while bounding
the number of choices by taking a factor of $\varDelta(\tvec)$ for
traversing an edge in $T$ and $\varDelta(S)$ for traversing an edge in $S-T$.
Finally, locations for the internal labelled vertices of $W^*$ can be chosen in
at most $w_{z,z}$ ways.

As we will show,
for each of conditions (2)--(4), we can choose $W^*$ and $\mathcal{Z}$ so 
that $\mathcal{Z}$ uses $\ell-z-2$ edges of~$T$ and $\ell-z-1$ edges of~$S-T$.

Condition (2) occurs if $\{v_i,u_{i+1}\}=\{j,k\}$ for some~$i$.
There are at most $\ell-1$ choices of~$i$, and two directions in which
$\{j,k\}$ can be traversed.
$W^*$ and $\mathcal{Z}$ can be chosen as follows.
\WZtable{
  $z+1\leq i$ & $W_1$ & $W[v_i,u_{z+1}], W[u_\ell,u_{i+1}]$ \\
  otherwise & $W_3$ & $W[u_{i+1},v_{\ell-z}], W[v_i,v_1]$}
Therefore,
\[
   N_2 \leq  2(\ell-1)  \varDelta(\tvec)^{\ell-z-2}\varDelta(S)^{\ell-z-1} w_{z,z} .
\]

Condition (3) occurs if $\{v_i,u_{i+1}\}=\{v_{i'},u_{i'+1}\}$ for some $i<i'$.
Note that $i'=i+1$ is impossible because condition (1) doesn't hold.
Allowing a factor of 2 for the two ways
that these steps can match, there are at most $(\ell-2)(\ell-3)$ choices.
$W^*$ and $\mathcal{Z}$ can be chosen as follows.
\WZtable{
  $z+1\leq i$ & $W_1$ & $W[v_\ell,v_{i+1}],  W[v_i,u_{z+1}]$ \\
   $2i+1\leq z+i+1\leq i'$ & $W_2(i)$ & $W[v_\ell,u_{i+z+1}], W[v_i,v_1]$  \\
  otherwise & $W_3$ & $W[u_1,u_{i'}], W[u_{i'+1},v_{\ell-z}]$}
Therefore
\[
    N_3 \leq (\ell-2)(\ell-3)  \varDelta(\tvec)^{\ell-z-2}\varDelta(S)^{\ell-z-1} w_{z,z}.
\]

Condition (4) occurs if $\{u_i,v_i\}=\{u_{i'},v_{i'}\}$ for some $i<i'$.
Considering that $i'=i+1$ is impossible and that there are two different orientations
for the overlap, the number of cases is bounded by~$(\ell-1)(\ell-2)$.
Also $(i,i')\ne (1,\ell)$, since we are assuming condition (1) doesn't hold.
$W^*$ and $\mathcal{Z}$ can be chosen as follows.
\WZtable{
  $z+1\leq i$ & $W_1$ & $W[v_\ell,u_{i+1}], W[u_i,u_{z+1}]$ \\
  $2i+1\leq z+i+1\leq i'$ and $i=1$ & $W_2(i)$ & 
     $W[u_1,u_{i+1}], W[u_{i'},u_{i+z+1}], W[v_\ell,u_{i'+1}]$ \\
  $2i+1\leq z+i+1\leq i'$ and $i\ne 1$ & $W_2(i)$ & 
     $W[v_\ell,u_{i+z+1}], W[u_1,v_{i-1}], W[v_i,u_{i+1}]$\\
  otherwise & $W_3$ &  $W[u_1,v_{i'-1}], W[v_{i'},v_{\ell-z}]$}
Therefore
\[
    N_4 \leq (\ell-1)(\ell-2) \varDelta(\tvec)^{\ell-z-1}\varDelta(S)^{\ell-z-2} w_{z,z}.
\]
Since $2(\ell-1)+(\ell-2)(\ell-1)+(\ell-2)(\ell-3)\leq 2\ell^2$,
part (a) of the lemma follows from~\eqref{e:pre-b}.

\medskip

Part (b) is proven in a completely similar way to  part (a). 
The number of walks  of length $2\ell-1$ 
from $j$ to $k$ which alternate between edges of $S$ and $T'$ starting from an edge in $T'$ equals  
$\langle(B'A)^{h-1}B' \evec_j , (B'A)^h B'\evec_k\rangle_S$
which is an upper bound for the number of  inverse $\ell$-switchings.
For the lower bound, we again need to consider four conditions when the constructed walk
$W$ is not a valid inverse $\ell$-switching: 
\begin{itemize}\itemsep=0pt

\item[($1'$)] One of the edges $v_i u_{i+1}$ which we choose from $S$ belongs also to  $T'$.
\item [($2'$)] Condition ($1'$) doesn't hold, but $\{u_i, v_i\} = \{j,k\}$ for some $i$. 
\item [($3'$)] Condition ($1'$) doesn't hold, but 
$\{u_i, v_i\} = \{u_{i'}, v_{i'}\}$ for some $i<i'$ (same $T'$ edge twice).
\item [($4'$)] Condition ($1'$) doesn't hold, but $\{v_i, u_{i+1}\} = \{v_{i'}, u_{i'+1}\}$
  for some $i<i'$ (same $S$ edge twice).

\end{itemize}

Let $N_1',\ldots,N_4'$ correspond to the counts for these cases.
Then we have 
\[ 
	N_1' \leq \ell w'_{\ell-2, \ell+1},
\]
because such walks consist of $\ell-2$ edges from $S$ and $\ell+1$ edges from $T'$,
and there are less than $\ell$ choices for which $S$ edge is in $T'$.
From now on we will assume that condition~($1'$) does not hold, so that walks
strictly alternate between $T'$ and $S-T'$. In particular, walks never
backtrack.
For each of the cases ($2'$)--($4'$), we will define $W^*$ and
$\mathcal{Z}$ such that $\mathcal{Z}$ uses $\ell-z-1$ edges of~$T'$ and
$\ell-z-2$ edges of~$S$.

For condition ($2'$), there are at most $2(\ell-1)$ possibilities for 
$i$ such that $\{v_i,u_{i+1}\}=\{j,k\}$.
$W^*$ and $\mathcal{Z}$ can be chosen as follows.
\WZtable{
  $z+1\leq i$ & $W_1$ & $W[v_i,u_{z+1}], W[u_{i+1},u_\ell]$ \\
  otherwise & $W_3$ & $W[v_i,v_1],W[u_{i+1},v_{\ell-z}]$}
Therefore,
\[
   N_2' \leq 2(\ell-1) \varDelta(\tvec)^{\ell-z-1}\varDelta(S)^{\ell-z-2} w'_{z,z}.
\]

For condition ($3'$), there are at most $(\ell-1)(\ell-2)$ possibilities
for $i<i'$ such that $i'\ne i+1$ and $\{u_i,v_i\}=\{u_{i'},v_{i'}\}$.
$W^*$ and $\mathcal{Z}$ can be chosen as follows.
\WZtable{
  $z+1\leq i$ & $W_1$ & $W[v_\ell,u_{i+1}], W[u_i,u_{z+1}]$ \\
  $2i+1\leq z+i+1\leq i'$ and $i=1$ & $W_2(i)$ & 
     $W[u_1,u_{i+1}], W[u_{i'},u_{i+z+1}], W[v_\ell,u_{i'+1}]$ \\
  $2i+1\leq z+i+1\leq i'$ and $i\ne1$ & $W_2(i)$ &
    $W[u_1,v_{i-1}], W[v_\ell,u_{i+z+1}], W[v_i,u_{i+1}]$ \\
  otherwise & $W_3$ & $W[u_1,v_{i'-1}], W[v_{i'},v_{\ell-z}]$}
Therefore,
\[
    N_3' \leq (\ell-1)(\ell-2) \varDelta(\tvec)^{\ell-z-1}\varDelta(S)^{\ell-z-2} w'_{z,z}.
\]

For condition ($4'$), there are at most $(\ell-2)(\ell-3)$ possibilities
for $i<i'$ such that $i'\ne i+1$ and $\{v_i,u_{i+1}\}=\{v_{i'},u_{i'+1}\}$.
$W^*$ and $\mathcal{Z}$ can be chosen as follows.
\WZtable{
  $z+1\leq i$ & $W_1$ & $W[v_\ell,v_{i+1}], W[v_i,u_{z+1}]$ \\
  $2i+1\leq z+i+1\leq i'$ & $W_2(i)$ & $W[u_1,u_{i+1}], W[v_{i'},u_{i+z+1}], W[u_\ell,u_{i'+1}]$ \\
  otherwise & $W_3$ & $W[u_1,u_{i'}], W[u_{i'+1},v_{\ell-z}]$}
Therefore,
\[
    N_4' \leq (\ell-2)(\ell-3) \varDelta(\tvec)^{\ell-z-1}\varDelta(S)^{\ell-z-2} w'_{z,z}.
\]
Since $2(\ell-1)+(\ell-1)(\ell-2)+(\ell-2)(\ell-3)\leq 2\ell^2$,
part (b) of the lemma follows.
This completes the proof.
\end{proof}

As a demonstration of the method, we start with the case of dense $S$ and then we proceed to 
Theorem \ref{l:co-sparse} for the sparse setting. 
	
	%\nicebreak
\subsection{Dense $S$}
 The following lemma estimates $\Pr(jk\in S_{\tvec})$ when $S$ is dense and $\tvec$ is a sparse degree sequence.
 \begin{lemma}\label{l:denseS}
		Let $\eps\in (0,1)$ be a constant
		and $S$ be a graph on $n$ vertices 
		such that 
		\[
		\varDelta(S) -\range{S} \geq \eps n,
		\qquad \text{and} \qquad 
		\left|\log \frac{\langle\xvec ,\yvec \rangle_S }{\|\xvec\|_1 \|\yvec\|_1 \varDelta(S)/n}\right| \leq \gamma,
		\]
		for all $\xvec,\yvec \in [0,1]^n$ with $\|\xvec\|_1, \|\yvec\|_1 \geq \eps^{7} n$. 
		Let  $\tvec$ be a degree sequence such that
		there exists a $\tvec$-factor of $S$ and
		\[
		\varDelta(\tvec) = o(n), \qquad  
		\frac{\varDelta(\tvec) - \range{\tvec}}{\varDelta(\tvec)} \geq \eps.
		\]
		Then, for any $jk \in S$, 
		\[
		\Pr(jk \in S_{\tvec}) = \exp\left(O
		\left( \gamma+ \dfrac{\range{S}}{\varDelta(S) } + \dfrac{\varDelta(\tvec)}{n} + \dfrac{\range{\tvec}}{\varDelta(\tvec)}
		\right)\right) \dfrac{\varDelta(\tvec)}{\varDelta(S)}.
		\]
  % \bdm{Would it be less disruptive to incorporate the new
  % error term by changing $\frac{\range{\tvec}}{\varDelta{\tvec}}$ into $\frac{\range{\tvec}+1}{\varDelta{\tvec}}$?}
  \end{lemma}
	
	\begin{proof}      
		Fix $\ell =7$. Using  Lemma \ref{l:FB}, we   will estimate 
		$|\mathcal{F}_\ell(T)|$ for $T\in \mathcal{S}(\tvec,jk)$ and 
		$|\mathcal{B}_\ell(T')|$ for $T'\in \mathcal{S}(\tvec,\overline{jk})$.
		Let $A=A(S)$, $B=A(T)$, $B'=A(T')$.     
		Let $s_{\min}$ and $t_{\min}$ denote $\varDelta(S)-\range{S}$ and $\varDelta(\tvec)-\range{\tvec}$ respectively.
		Using the assumptions of the lemma and $\|A \evec_j\|_\infty\leq 1$, and noting that
		$\|B\|_{\infty}\leq \varDelta(\tvec)$ and $\|A\|_{\infty}\leq \varDelta(S)$, we find, for $h=2,3$,
		\[
		\frac{\|(BA)^h \evec_j\|_1}
		{\|(BA)^h \evec_j\|_\infty} \geq \frac{(t_{\min} s_{\min})^h}{
			\|B\|_\infty^h \|A\|_{\infty}^{h-1} \|A \evec_j\|_\infty} \geq \eps^{2h-1} s_{\min} \geq  \eps^{6}n.
		\]
    The first inequality above holds since $s_{\min}^h t_{\min}^h$ is a lower bound on the number of walks of $h$ steps in $T$ and $h$ steps in $S$, starting at vertex~$j$, whereas $$\|(BA)^h \evec_j\|_\infty\leq \norm{(BA)^{h-1}B}_{\infty}\norm{A\evec_j}_{\infty} \leq \|B\|_\infty^h \|A\|_{\infty}^{h-1} \|A \evec_j\|_\infty.$$
    Using $\norm{A\evec_j}_{\infty}\leq 1$, and writing $\vvec=(v_1,\ldots,v_n)=\sum_{i=1}^n v_i \evec_i$ for any $\vvec\in {\mathbb R}^n$,
    \[
    \norm{A\vvec}_{\infty}\leq \sum_{i=1}^n |v_i|\norm{A\evec_i}_{\infty} \leq \sum_{i=1}^n |v_i|=\norm{\vvec}_1.
    \]
    Thus, we have 
		\[
		\frac{\|(B'A)^h B' \evec_j\|_1}
		{\|(BA)^h  B'\evec_j\|_\infty}  \geq 
		\frac{(t_{\min} s_{\min})^h  t_{\min}}{
			\|B'\|_\infty^h \|A\|_{\infty}^{h-1} \|AB' \evec_j\|_\infty}
		\geq \eps^{2h-1} \frac{s_{\min} t_{\min}}{\|B' \evec_j\|_1} 
		\geq \eps^{2h} s_{\min} \geq \eps^{7} n.
		\]
		Note also that 
		\begin{align*}
			\|(BA)^h \evec_j\|_1 &=  (\varDelta(\tvec) \varDelta(S))^h  
			\exp\left( O\left(
			\dfrac{\range{S}}{s_{\min}} + \dfrac{\range{\tvec}}{t_{\min}}
			\right)\right),\\
			\|(B'A)^h B' \evec_j\|_1 &=  (\varDelta(\tvec))^{h+1} (\varDelta(S))^h  
			\exp\left( O\left(
			\dfrac{\range{S}}{s_{\min}} + \dfrac{\range{\tvec}}{t_{\min}}
			\right)\right).
		\end{align*}
		Using similar bounds for $\evec_k$ and the assumption of the lemma on $\langle \xvec, \yvec \rangle_S$, we find that 
		\begin{align*}
			\langle(BA)^{3} \evec_j , 
			(BA)^{3} \evec_k\rangle_S &= 
			\exp\left( O\left( \gamma+
			\dfrac{\range{S}}{s_{\min}} + \dfrac{\range{\tvec}}{t_{\min}}
			\right)\right) \varDelta(\tvec)^{6} \varDelta(S)^{7}/n,
			\\
			\langle(B'A)^{2}B' \evec_j , 
			(B'A)^3 B'\evec_k\rangle_S 
			&= 
			\exp\left( O\left( \gamma+
			\dfrac{\range{S}}{s_{\min}} + \dfrac{\range{\tvec}}{t_{\min}}
			\right)\right) \varDelta(\tvec)^{7} \varDelta(S)^{6}/n.
		\end{align*}     
		Next, we need to bound the quantities $w_{a,b}$ and $w_{a,b}'$ that appear in the 
		lower bounds of Lemma~\ref{l:FB}. 
		Consider any product $P$ of $a\geq 1$ factors $A$ and $b$ factors $B$. 
		Representing $P= P_1 A P_2$, we get that 
		\[
		\|P_1 A P_2 \evec_i\|_\infty \leq \|P_1\|_\infty \|P_2 \evec_i\|_1 \leq \varDelta(S)^{a-1} \varDelta(\tvec)^{b} \leq    \varDelta(S)^a \varDelta(\tvec)^b/\eps n.
		\]
		Thus, we can bound  $w_{6,7} = O( \varDelta(\tvec)^7\varDelta(S)^6/n)$, 
		$w_{5,8}' = O(\varDelta(\tvec)^8\varDelta(S)^5 /n) $, and
		$w_{2,2}, w'_{2,2} = O(\varDelta(\tvec)^2\varDelta(S)^2/n)$.
		Thus, applying Lemma \ref{l:FB}, we conclude that 
		\begin{align*}
			|\mathcal{F}_\ell(T)| &= 
			\exp\left( O\left( \gamma  + 
			\dfrac{\range{S}}{s_{\min}} + \dfrac{\range{\tvec}}{t_{\min}}
			+ \dfrac{\varDelta(\tvec)}{n}\right)\right) \varDelta(\tvec)^6 \varDelta(S)^7/n,\\
			|\mathcal{B}_\ell(T')| &=
			\exp\left( O\left(\gamma  + 
			\dfrac{\range{S}}{s_{\min}} + \dfrac{\range{\tvec}}{t_{\min}}
			+ \dfrac{\varDelta(\tvec)}{n}
           \right)\right) \varDelta(\tvec)^7 \varDelta(S)^6/n.
		\end{align*}
		Combining \eqref{switch1} and \eqref{switch2} completes the proof.
	\end{proof}

	%\nicebreak
	\subsection{Preliminaries for sparse $S$}
	%%%%%%%%%%%%%%%%%%%%%%%%%%%%%%

	For a sparse  $S$,  estimating $|\mathcal{F}_\ell|$ and $|\mathcal{B}_\ell|$ accurately is a non-trivial task and relies heavily on the pseudorandom properties of $S$. Here we prove bounds 
	for the quantities $w_{a,b}(S,T)$ which appear in Lemma \ref{l:FB}. 
	First we consider the case when both graphs are regular.
	Let $J$ denote the $n\times n$ matrix where every entry is 1.

	\begin{lemma}\label{lem:walk}
		Let $T$ be a $t$-regular graph and  $S$ be a  $s$-regular graph on the same vertex set~$[n]$.
		Assume that $\bigl\|A({S}) - \frac{|S|}{N} J\bigr\|_2 \leq {s} n^{-\alpha}$ for some $\alpha>0$.
		Then, for any integers $a\geq 1/\alpha$, $b\geq 0$, we have
		\[
		w_{a,b}({S},{T})\leq \frac{2s^at^b}{n}.
		\]
	\end{lemma}
	
	\begin{proof}
		Let $A=A(S)$ and $B = A(T)$.
		Consider any matrix $P$ which is a product of $a$ factors $A$ and $b$ factors $B$.
		Let  $\tilde{P}$ denote the matrix obtained by replacing all factors $A$ in $P$ by   $A -  \frac{{|S|}}{N} J$. 
		Write 	$\evec_i = \boldsymbol{1}/n + \vvec$, where $\boldsymbol{1}$ is the vector with all components equal $1$.  Note that 
		$A \boldsymbol{1} = {s}  \boldsymbol{1}$,  
		$B \boldsymbol{1} = {t}  \boldsymbol{1}$
		and  $\vvec \perp  \boldsymbol{1}$.  
		Since  operators $A$ and $A -  \frac{|S|}{N} J$ act identically on the space orthogonal to $\boldsymbol{1}$,  we find  that
		$
		P\vvec = \tilde{P}\vvec.
		$
  Since $B$ is symmetric $\norm{B}_2$ is equal to the largest eigenvalue of $B$, which is equal to $t$. 
		Thus, we obtain that
		\begin{align}
			\| P \evec_{i}\|_\infty 
			&\leq \dfrac{\norm{P\boldsymbol 1}_\infty}{n} + 
			\| P \vvec\|_\infty=  \dfrac{\norm{P\boldsymbol 1}_\infty}{n} + 
			\| \tilde{P} \vvec\|_\infty \leq \dfrac{\norm{P\boldsymbol{1}}_\infty}{n} + 
			\| \tilde{P} \vvec\|_2 \nonumber \\
			& \leq \dfrac{s^a t^b}{n} +  \| (A- \tfrac{|S|}{N}J)\|_2^a \cdot \|B\|_2^b 		\leq  \dfrac{s^a t^b}{n} + s^a t^b n^{-a \alpha} \leq \dfrac{2s^a t^b}{n}.\label{eq:inequality-walk-counts} 
		\end{align}
Again, the first inequality in~\eqref{eq:inequality-walk-counts} above holds since $\norm{P\boldsymbol 1}_{\infty}$ denotes the maximum number of walks starting from a given vertex, which use exactly $a$ edges in $S$ and $b$ edges in $T$, and since $\norm{\vvec}_2^2=(1-1/n)^2+(n-1)/n^2=(n-1)/n<1$,
  we have that $\norm{\tilde{P} \vvec}_2\leq \norm{\tilde{P}}_2 \leq \| (A- \tfrac{|S|}{N}J)\|_2^a \cdot \|B\|_2^b $.
		Taking the maximum over all $i$ and $P$ completes the proof.
	\end{proof}

	We will need a bound similar to Lemma \ref{lem:walk} for non-regular $S$ and $T$ as well. For this purpose, we construct regular supergraphs $\tilde{S}\supseteq S$ and $\tilde{T} \supseteq T$ and estimate 
	\begin{equation}\label{eq:wSS}
		w_{a,b}(S,T) \leq w_{a,b}(\tilde{S},\tilde{T}) .
	\end{equation}
	The next lemma shows that if $G$ is a graph with small $\range{G}$  then there exists a regular supergraph $\tilde{G}\supseteq G$ which is not much bigger than $G$.
	\begin{lemma}\label{l:supergraph}
		Let $G$ be a graph on $n$ vertices such that $\varDelta(G)+3\range{G} < n/4$.
		If $d$ is an even number that
		\[ 
		\varDelta(G) + \range{G}  \leq d \leq \varDelta(G)+2\range{G}
		\] 
		then there exists a $d$-regular supergraph $\tilde{G}$  of $G$. 
	\end{lemma}
	\begin{proof}
		If $\range{G}=0$ there is nothing to prove as we can take $\tilde{G} = G$. Thus, we may assume otherwise. Define sequence $\rvec = (r_1,\ldots,r_n)$ by $r_i=d-d_G(i)$, where $d_G(i)$ is the degree of vertex $i$ in $G$. 
		It is sufficient to  find an $\rvec$-factor of $K_n-G$ because  the union of this $\rvec$-factor and $G$ gives our desired $\tilde{G}$. 
		
		By the assumptions, for all vertices $i$, we have
		\begin{equation}\label{rran}
			\range{G}\leq d - \varDelta(G)  \leq  r_i \leq d - \varDelta(G) + \range{G} \leq 3 \range{G}. 
		\end{equation}
		
		Then, for any $U \subseteq [n]$,   we have
		\[	
		\sum_{i \in U} r_i \leq 3 \range{G}\cdot |U|  \leq 
		|U|\cdot(|U|-1) 
		+ \sum_{i \notin U} \min\{|U|, r_i\}.
		\]
		To see the above inequality holds, note that if $|U| \geq 3 \range{G}+1$ then the first term of the RHS is at least $3 \range{G}\, |U| $.
		If, on the other hand, $|U|\leq 3\range{G}< n/4$ then the second term of the RHS
		is at least $\frac{3}{4}n\min\{|U|,\range{G}\}\geq 3 \range{G}\, |U|$.
		Also, $\sum_{i} r_i$ 
		is even since $d$ is even. By the Erd\H{o}s--Gallai theorem we conclude that 
		${\rvec}$ is a graphical degree sequence.

		Let $R$ be an $\rvec$-factor of $K_n$ such that $R \cap G$ has the smallest number of edges. We use a switching-type argument to show that that $R \subseteq K_n - G$. By contradiction, assume that there is an edge $u_1 v_1 \in R \cap G$.  
		Consider  edges $u_2 v_2 \in R$ 
		such that $u_1 u_2 \in K_n -  (R\cup G)$.
		The number of choices for such $u_2 v_2$ is at least
		\[	
		\bigl(n-2- (\varDelta(G)-1) - (\varDelta(R)-1)\bigr)  \range{G}
		= (n- \varDelta(G) - \varDelta(R))  \range{G}.
		\] 
		The first factor in the LHS  corresponds to choices $u_2$ that  $u_1 u_2 \notin R\cup G$ and  the second factor  in the LHS is  a lower bound for the number of ways to choose $v_2$ given $u_2$ (by~\eqref{rran}).  
		Note that among all the choices for $u_2v_2$ above, at most $(\varDelta(G) + \varDelta(R)-2)\varDelta(R)$ choices satisfy $v_1 v_2 \in R\cup G$ (estimating the number of ways to choose $v_2$ and then $u_2$).
		Also, among all the choices for $u_2v_2$ above, at most $\varDelta(R)$ 
		choices satisfy $v_1 = v_2$.
		By the assumptions and \eqref{rran},
		we find that
		\begin{align*}
			(n- \varDelta(G) - \varDelta(R))  \range{G} >  (\varDelta(G) + \varDelta(R)-1)\varDelta(R).
		\end{align*}
		Therefore, we can find such $u_2 v_2 \in R$ that  $u_1 u_2 \notin R\cup G $ and 
		$v_1 v_2 \notin R\cup G$ and all vertices $u_1,v_1,u_2,v_2$ are distinct.  
		Then we can replace edges $u_1 v_1$ and $u_2 v_2$ by $u_1 u_2$ and $v_1 v_2$ to get an 
		$\rvec$-factor which has fewer common edges with $G$ than $R$ does. This contradicts our choice of  $R$.  Therefore $R \cap G$ must be empty, which completes the proof.
	\end{proof}

	%%%%%%%%%%%%%%%%%%%%%%%%%%%%%%%%%
	%%%%%%%%%%%%%%%%%%%%%%%%%%%%%%%%%
	
	\nicebreak
	\subsection{Proof of Theorem \ref{l:co-sparse}.}
 \label{sec:proof_l_co-sparse}

	First we  consider the case when $\varDelta(S) \geq n/16$.  In this case, we show that the required probability bound \eqref{eq:edge-prob}   holds with 
      $\delta = \eps^7$, where  $\eps = \min\{1/32,1/\beta\}$.  We can assume that $\gamma_S(\delta)<\infty$ since 
      \eqref{eq:edge-prob} is trivial otherwise.
    The proof is by applying Lemma~\ref{l:denseS}.    
  Indeed, from (A1)  we get that $\varDelta(\tvec)=o(n)$.
  Assumptions (A2) and (A3) of Theorem~\ref{l:co-sparse}  imply that $\range{S} \ll \varDelta(S)$: if $\alpha = \Omega(1)$ then it follows from (A2) and if $\alpha \rightarrow 0$ then it follows (A3).
  Thus,  we get 
  \[
    \varDelta(S)-\range{S} \geq n/16 - o(n)\geq  \eps n
  \]
for any $\eps \leq 1/32$.  Also, from (A3)  we get that
$
    \dfrac{\varDelta(\tvec)} {\varDelta(\tvec)-\range{\tvec}} \leq \beta.
$ 
Thus all the assumptions of  Lemma~\ref{l:denseS} hold
and  applying this lemma gives the required bound~\eqref{eq:edge-prob}.

\smallskip
	In the following, we assume that 
	$\varDelta(S) < n/16$ which implies 
    $
    \varDelta(S) + 3\range{S}  <n/4
    $
	and  \[\varDelta(\tvec) + 3\range{\tvec}<
     4 \varDelta(\tvec) \leq  4 \varDelta(S)<n/4.
     \]
     Take $\ell$ to be the odd number  
	in $\{\lceil 3/\alpha\rceil+3,\lceil 3/\alpha\rceil+4\}$. Using  Lemma \ref{l:FB}, we  will estimate 
	$|\mathcal{F}_\ell(T)|$ for $T\in \mathcal{S}(\tvec,jk)$ and 
	$|\mathcal{B}_\ell(T')|$ for $T'\in \mathcal{S}(\tvec,\overline{jk})$.
	Let $A=A(S)$, $B=A(T)$, $B'=A(T')$.        
	Recall that $s_{\min}=  \varDelta(S) - \range{S}$ is the smallest degree of $S$      and $t_{\min} = \varDelta(\tvec) -\range{\tvec}$   is the smallest component of $\tvec$. 
Note that assumption (A3) implies that $t_{\min}>0$.
	Let $\tilde{S}$, $\tilde{T}$  be regular supergraphs of $S$ and $T$ given by Lemma~\ref{l:supergraph}. We have 
	\begin{equation}
	\frac{ \varDelta(\tilde{S})}
	{s_{\min}} \leq  
	\frac{\varDelta(S) + 2 \range{S}}{s_{\min}}
	= 1  + \frac{3 (\varDelta(S) - s_{\min}) }{s_{\min}} \leq \left( \frac{\varDelta(S)}{s_{\min}}\right)^3.\label{eq:ratioStilde}
	\end{equation}
	Similarly, $\varDelta(\tilde{T})/t_{\min} \leq (\varDelta(\tvec)/t_{\min})^3$.   
	Applying bound \eqref{eq:Gspec} for $G =\tilde{S}\setminus S$ and
	recalling from Lemma \ref{l:supergraph} that $\tilde{S}$ is regular and	$\varDelta(\tilde{S}) \leq \varDelta(S) + 2\range{S},$
	we find that
	\begin{align*}
		\|A(\tilde{S}) - \dfrac{|\tilde{S}|}{N} J\|_2
		&\leq \|A(S) -  \dfrac{|S|}{N} J\|_2 +  \|A(\tilde{S}\setminus S) - \dfrac{|\tilde S\setminus S|}{N} J\|_2 \\
		&\leq  \|A(S) -  \dfrac{|S|}{N} J\|_2  + 2   \varDelta(\tilde{S}\setminus S)\\
		&\leq   \|A(S) - \dfrac{|S|}{N} J\|_2 +   2 (\varDelta(\tilde{S}) - s_{\min}) 
		\\		
		&\leq   \|A(S) - \dfrac{|S|}{N} J\|_2 +   6 \range{S}  \\ &\leq n^{-\alpha}\varDelta(S)\leq n^{-\alpha}\varDelta(\tilde S),\quad \mbox{by (A2)}.
	\end{align*}
	Combining Lemma \ref{lem:walk} and estimate \eqref{eq:wSS}, we find that, for $a \geq 1/\alpha$, 
	\begin{equation}\label{bounds_wab}
		\begin{aligned}
			w_{a,b}(S,T)&\leq  w_{a,b}(\tilde{S},\tilde{T}) \leq \frac{2 \varDelta(\tilde{S})^a  \varDelta(\tilde{T})^b}
			{n}
			\\ 
			&\leq \frac{2 s_{\min}^a t_{\min}^b}{n} 
			\left( \frac{\varDelta(\tilde{S})}{s_{\min}} \cdot \frac{\varDelta(\tilde T)}{t_{\min}} \right)^{\max\{a,b\}} 
			\leq \frac{2  s_{\min}^a t_{\min}^b}{n} \beta^{3 \max\{a,b\} \alpha},
		\end{aligned}
	\end{equation}
 where the last inequality above holds by~\eqref{eq:ratioStilde} and  (A3).
	The same bound holds for $w_{a,b}(S,T')$.
Observing   that $\ell=O(1/\alpha)$, we get that, for $1/\alpha \leq h \leq (\ell-1)/2$,  
	\begin{align*}
		\frac{\|(BA)^h \evec_j\|_1}
		{\|(BA)^h \evec_j\|_\infty} &\geq 
		\frac{(s_{\min} t_{\min})^h}{w_{h,h}} \geq 
		n \beta^{ -3 h \alpha}/2 \geq \delta n,\\
		\frac{\|(B'A)^h B' \evec_j\|_1}
		{\|(BA)^h  B'\evec_j\|_\infty}  &\geq 
		\frac{s_{\min}^h  t_{\min}^{h+1}}{w_{h,h}} \geq 
		n \beta^{ -3 (h+1) \alpha}/2 
        \geq \delta n,
	\end{align*}
    where $\delta = \delta(\beta)>0$ corresponds to the largest $h$ within the specified range since $\beta\geq 1$. We also get
	\begin{align*}
		\|(BA)^h \evec_j\|_1 &=  (\varDelta(\tvec) \varDelta(S))^h  
		\exp\left( O\left(
		\dfrac{\range{S}}{\alpha s_{\min}} + \dfrac{\range{\tvec}}{\alpha t_{\min}}
		\right)\right),\\
		\|(B'A)^h B' \evec_j\|_1 &=  (\varDelta(\tvec))^{h+1} (\varDelta(S))^h  
		\exp\left( O\left(
		\dfrac{\range{S}}{\alpha s_{\min}} + \dfrac{\range{\tvec}}{\alpha t_{\min}}
		\right)\right).
	\end{align*}
	Using   similar bounds for $\evec_k$ and recalling the definition of $\gamma_S$ from \eqref{def:gamma-S} (and  scaling by the infinity norm), we find that
	\begin{align*}
		\langle(BA)^{\frac{\ell-1}{2}} \evec_j , 
		(BA)^{\frac{\ell-1}{2}} \evec_k\rangle_S &= 
		\exp\left( O\left(\gamma_S(\delta) + 
		\dfrac{\range{S}}{\alpha s_{\min}} + \dfrac{\range{\tvec}}{\alpha t_{\min}}
		\right)\right) \varDelta(\tvec)^{\ell-1} \varDelta(S)^{\ell}/n,
		\\
		\langle(B'A)^{\frac{\ell-3}{2}} B' \evec_j , 
		(B'A)^{\frac{\ell-1}{2}} B'\evec_k\rangle_S 
		&= 
		\exp\left( O\left( \gamma_S(\delta) + 
		\dfrac{\range{S}}{\alpha  s_{\min}} + \dfrac{\range{\tvec}}{\alpha t_{\min}}
		\right)\right) \varDelta(\tvec)^{\ell} \varDelta(S)^{\ell-1}/n.
	\end{align*}     
	Next, using again $\ell=O(\alpha^{-1})$, we apply \eqref{bounds_wab} to bound 
	quantities $w_{a,b}$ and $w_{a,b}'$ that appear in the 
	lower bounds of Lemma \ref{l:FB}. Thus, we find that 
		\begin{align*}
		\ell w_{\ell-1,\ell} &= O\Bigl(\dfrac{ \varDelta(\tvec)^{\ell}\varDelta(S)^{\ell-1}}{\alpha n}\Bigr),
		 & \ell^2 \varDelta(\tvec)^{\ell-z-2}\varDelta(S)^{\ell-z-1}w_{z,z}
		   = O\Bigl( \dfrac{\varDelta(\tvec)^{\ell-2}\varDelta(S)^{\ell-1}}{\alpha^2 n}\Bigr),
		\\
		\ell w_{\ell-2,\ell+1} &=  O\Bigl(\dfrac{\varDelta(\tvec)^{\ell+1}\varDelta(S)^{\ell-2} }{\alpha n}\Bigr),
		& \ell^2 \varDelta(\tvec)^{\ell-z-1}\varDelta(S)^{\ell-z-2}w'_{z,z}
		   = O\Bigl( \dfrac{\varDelta(\tvec)^{\ell-1}\varDelta(S)^{\ell-2}}{\alpha^2 n}\Bigr).
	\end{align*}
	Applying Lemma \ref{l:FB}, we conclude that 
	\begin{align*}
		\frac{|\mathcal{B}_\ell(T)|}{|\mathcal{F}_\ell(T')| } = 
		\exp\left( O\left(\gamma_S(\delta)  +
		\dfrac{\range{S}}{\alpha  s_{\min}} + \dfrac{\range{\tvec}}{\alpha t_{\min}}
		+ \dfrac{\varDelta(\tvec)}{\alpha \varDelta(S)}
		+ \dfrac{1}{\alpha^2\varDelta(\tvec)\varDelta(S)}\right) \right) 
		\frac{\varDelta(\tvec)}{\varDelta(S)}.
	\end{align*}
	Combining  \eqref{switch1} and \eqref{switch2} completes the proof of~\eqref{eq:edge-prob}.

	%%%%%%%%%%%%%%%%%%%%%%%%%%%%%%%%%
	%%%%%%%%%%%%%%%%%%%%%%%%%%%%%%%%%

	\section{Appendix}\label{sec:appendix}
	
	Here we prove or cite the technical lemmas that are used in the proofs.
	This section  is self-contained and does not rely on assumptions other than those stated.

	\begin{lemma}\label{l:technical}
		Let  $X$ be a binomial, Poisson, or hypergeometric random variable. Then,
		for any $\eps\geq 0$,
		\[
		\Pr( |X -\E X|  \geq  \eps \E X  ) \leq 2  e^{-\frac{\eps^2}{2+\eps} \E X }.
		\]
	\end{lemma}
	\begin{proof}
	  This follows from~\cite[Theorem 2.1, Remark 2.6,  and Theorem 2.10]{JansonBook}. 
	\end{proof}

	\begin{lemma}\label{A:spec}
		Let $A_p$ be the adjacency matrix of $G\sim G(n,p)$, where
		$\log n\leq pn\leq n-\log n$.  Assume 
		$\eps = \eps(n)$ satisfies
		$ \sqrt{\dfrac{\log n}{pn}}\ll \eps \ll 1$. Then 
		\[
		\Pr\left(\|A_p - pJ \|_2  \leq \eps p n\right)
		=1 - e^{-\Omega(\eps ^2 p^2 n^2)},
		\]
		where $J$ denotes the $n\times n$ matrix with all entries equal $1$.
	\end{lemma}

\begin{proof}
Let $\lambda(A_p)=\norm{A_p-pJ}_2$. This is the spectral norm of $A_p$ since the matrix is symmetric.
Since $\E A_p=pJ-pI$ and $p\ll\eps pn$ under our assumptions, we have $\lambda(A_p)
\leq \norm{A_p-\E A_p}_2+\norm{pI}_2 
= \norm{A_p-\E A_p}_2+o(\eps pn)$.

First consider $\log^4 n\leq pn\leq 1-\log^4 n$.
By~\cite[Theorem 1.2]{Vu2007}, $\lambda(A_p) \leq \E \lambda(A_p)
+ \frac12\eps pn$ with probability at least $1 - e^{-\Omega(\eps ^2 p^2 n^2)}$.
Also, under our assumption on $p$, \cite[Theorem 1.4]{Vu2007}
gives $\E \lambda(A_p)\leq 2p^{1/2}(1-p)^{1/2}\sqrt{n}+Cp^{1/4}n^{1/4}\log n+p$ for some constant~$C$.
Using our lower bound on $\eps$,
this implies that $\E\lambda(A_p)=o(\eps pn)$.
Consequently, $\lambda(A_p)\leq\eps pn$ with probability
at least $1 - e^{-\Omega(\eps ^2 p^2 n^2)}$.

For $\log n\leq pn\leq\log^4 n$, we can use~\cite[Corollary 3.3]{BGBK}.
Note that the parameter~$d$ in~\cite{BGBK} equals $(n-1)p$ and that
our lower bound on $\eps$ implies $\eps^{-2}\log^{1/2}n \,(\log\log n)^{-1/2}=o(pn)$,
so~\cite[Eq. (3.7)]{BGBK} is satisfied.
Therefore, we have $\lambda(A)\leq \frac12\eps pn+p\leq \eps pn$ with probability at least $1-e^{-\Omega(\eps^2p^3n^3)}
=1-e^{-\Omega(\eps^2p^2n^2)}$.

For $n-\log^4 n\leq pn\leq n-\log n$, note that
$A_p-\E A_p=-(A_{1-p}-\E A_{1-p})$.
Applying~\cite[Corollary 3.3]{BGBK} to $A_{1-p}$, we have
$\Pr\bigl(\norm{A_p-\E A_p}_2>\frac12\eps pn\bigr)
= e^{-\Omega((1-p)\eps^2n^3)}
= e^{-\Omega(\eps^2p^2n^2)}$ since $1-p\gg n^{-1}$
and $p\sim1$.
Since $\eps pn\gg 1$, the same bound holds
for $\Pr\bigl(\norm{A_p-pJ}_2>\eps pn\bigr)$.
\end{proof}

	Next, we prove that binomial random graphs are pseudorandom 
	in the following sense. Recall from \eqref{Z_def} that, for $\xvec,\yvec \in \Reals^n$ and a graph $G$, 
	\[		\langle  \xvec,\yvec\rangle_G = \sum_{(jk)\st jk \in G} x_j y_k.
	\]
	\begin{lemma}\label{A:z_con}
		Let  $G \sim \G(n,p)$ for some $1\geq p \gg \log n/n$.  
		Assume 
		$\eps = \eps(n)>0$  such that
		$ \sqrt{\dfrac{\log n}{pn}}\ll \eps \ll 1$.
		Then, with probability $1- e^{-\Omega(\eps^2 p n^2)}$,  we have 
		\[
		\bigl|\langle \xvec, \yvec \rangle_G  -  p \|\xvec\|_1  \|\yvec\|_1 \bigr| \leq  \eps p \|\xvec\|_1  \|\yvec\|_1
		\] 
		uniformly for all  $\xvec, \yvec \in [0,1]^n$  with $\|\xvec\|_1,\|\yvec\|_1 = \Omega(n)$.
	\end{lemma}
	\begin{proof}
		For any
		$\xvec, \yvec \in [0,1]^n$  with $\|\xvec\|_1,\|\yvec\|_1  = \Omega(n)$,  we have that
		\[
		\E \langle \xvec, \yvec \rangle_G= p \|\xvec\|_1 \|\yvec\|_1 -  p \sum_{j=1}^n x_j y_j = 
		\left(1+ O(n^{-1})\right)p \|\xvec\|_1 \|\yvec\|_1.
		\]
		Note that $\Var \langle \xvec, \yvec \rangle_G  \leq 2pn^2$.
		Using McDiarmid's inequality \cite[Theorem 2.7]{McDiarmid} (with $V=\Var \langle \xvec, \yvec \rangle_G$, $b=2$ and $t=\eps p \|\xvec\|_1  \|\yvec\|_1/2$), we  get
		\begin{align*}
			\Pr \left( \frac{|\langle \xvec, \yvec \rangle_G- \E \langle \xvec, \yvec \rangle_G|}{p \|\xvec\|_1  \|\yvec\|_1} \geq \eps /2  \right)
			&\leq  2\exp \left( - \frac{ (\eps p \|\xvec\|_1  \|\yvec\|_1/2)^2 }{ 2 \Var \langle \xvec, \yvec \rangle_G+ (2/3)  \eps p \|\xvec\|_1  \|\yvec\|_1  }\right)
			\\ &= 		e^{-\Omega(\eps^2 p n^2)}.
		\end{align*}	
		To make the probability estimate hold for all such $\xvec$, $\yvec$, we approximate them 
		with $\xvec',\yvec' \in \{ j/n \st j=1,\ldots,n\}^n$ such that 
		$\|\xvec - \xvec'\|_\infty \leq n^{-1}$ and $\|\yvec - \yvec'\|_\infty \leq n^{-1}$. 
		Denoting $\boldsymbol{1} = (1,\ldots,1)\trans $, we find that   

\begin{align*}
		|\langle \xvec, \yvec \rangle_G- \langle \xvec', \yvec' \rangle_G |
		&\leq |\langle \xvec, \yvec \rangle_G- \langle \xvec, \yvec' \rangle_G |+|\langle \xvec, \yvec' \rangle_G- \langle \xvec', \yvec' \rangle_G |\\
  &=|\langle \xvec, \yvec-\yvec' \rangle_G |+|\langle \xvec-\xvec', \yvec' \rangle_G | \\ &=O(n^{-1})\left(\langle \xvec,\boldsymbol{1}\rangle_G + \langle \yvec',\boldsymbol{1}\rangle_G \right)=O(n),
		\end{align*}
  as
		$ \langle\xvec , \boldsymbol{1}\rangle_G , \langle \yvec' , \boldsymbol{1}, \rangle_G \leq \langle \boldsymbol{1}, \boldsymbol{1}\rangle_G  \leq n^2$. Observe that $\eps p \|\xvec\|_1  \|\yvec\|_1  \gg n$, so all the error terms are within the required range.
		Allowing $n^{2n}$ for choice of $\xvec',\yvec'$, using the union bound
		and recalling that $n \log n \ll \eps^2 p n^2$, the proof is complete.
	\end{proof}	
	
	%%%%%%%%%%%%%%%%%%%%%%%%%%%%%%%%%
	%%%%%%%%%%%%%%%%%%%%%%%%%%%%%%%%%

\section*{Acknowledgment} 

We especially wish to thank the referees who devoted an
exceptional amount of effort to this paper and assisted
us considerably towards the final version.
		
	\nicebreak
	 
	%%%%%%%%%%%%%%%%%%%%%%%%%%%%%%%%%%%%%%%%%%%%%%%%%%%%%%%%%%%%%%%
	
\end{document}